\documentclass[lettersize,journal]{IEEEtran}

\usepackage{amsmath,amssymb,amsfonts}
\usepackage{mleftright}
\usepackage{ifthen}   
\usepackage{graphicx} 
\usepackage[mode=buildnew]{standalone}  
\usepackage{pgfplots} 
\pgfplotsset{table/search path={inc},compat=1.16}
\usepackage{tikzscale}
\usepackage{balance}
\usepackage[hidelinks]{hyperref} 
\usepackage{booktabs} 
\usepackage{makecell}
\usepackage{multirow}
\usepackage{subcaption}  
\usepackage{enumitem} 
\usepackage{float}

\newcommand{\bma}{\begin{bmatrix}}
\newcommand{\ema}{\end{bmatrix}}

\newcommand{\diag}[1]{\operatorname{diag}\left(#1\right)} 
\newcommand{\T}{{\mathsf{T}}} 

\newcommand{\Reals}{\mathbb{R}}      

\newcommand{\Normal}[1]{\mathcal{N}\!\left({#1}\right)} 

\newcommand{\dd}{\mathop{}\!\mathrm{d}} 

\newcommand{\eqdef}{\triangleq} 

\newcommand{\argmin}{\operatorname*{argmin}}
\newcommand{\argmax}{\operatorname*{argmax}}

\makeatletter
\DeclareFontFamily{U}{MnSymbolA}{}
\DeclareSymbolFont{MnSyA}{U}{MnSymbolA}{m}{n}
\DeclareFontShape{U}{MnSymbolA}{m}{n}{
<-6> MnSymbolA5
<6-7> MnSymbolA6
<7-8> MnSymbolA7
<8-9> MnSymbolA8
<9-10> MnSymbolA9
<10-12> MnSymbolA10
<12-> MnSymbolA12}{}
\DeclareMathSymbol{\smallrightarrow}{\mathrel}{MnSyA}{0}
\DeclareMathSymbol{\smallleftarrow}{\mathrel}{MnSyA}{2}
\DeclareMathSymbol{\smallleftrightarrow}{\mathrel}{MnSyA}{16}
\newcommand{\smallrightarrowfill@}{\arrowfill@\relbar\relbar\smallrightarrow}
\newcommand{\smallleftarrowfill@}{\arrowfill@\smallleftarrow\relbar\relbar}
\newcommand{\smallleftrightarrowfill@}
{\arrowfill@\smallleftarrow\relbar\smallrightarrow}
\renewcommand{\overrightarrow}{\mathpalette{\overarrow@\smallrightarrowfill@}}
\renewcommand{\overleftarrow}{\mathpalette{\overarrow@\smallleftarrowfill@}}
\renewcommand{\overleftrightarrow}
{\mathpalette{\overarrow@\smallleftrightarrowfill@}}
\makeatother
\providecommand{\msgf}[2]{\protect\overrightarrow{#1}_{\mspace{-3mu}#2}} 
\providecommand{\msgb}[2]{\protect\overleftarrow{#1}_{\mspace{-3mu}#2}} 

\DeclareMathOperator{\E}{\textnormal{\ensuremath{\mathbb{E}}}}
\newcommand{\EE}[1]{\E\!\left[{#1}\right]}

\newcommand{\Var}[1]{\operatorname{Var}\left[#1\right]} 

\newcommand{\cond}{\hspace{0.02em}|\hspace{0.08em}}

\newcounter{examplecntr}
\newenvironment{example}[1][]%
{\begin{trivlist}\small\item[]\refstepcounter{examplecntr}%
 {\bfseries Example~\theexamplecntr%
  \ifthenelse{\equal{#1}{}}{}{ (#1)}.
}}%
{\end{trivlist}}

\newcounter{definitioncntr}
{\begin{trivlist}\item[]\refstepcounter{definitioncntr}%
{\bfseries Definition~\thedefinitioncntr.}}%
{\hfill$\Box$\end{trivlist}}

\newcounter{theoremcntr}
\newenvironment{theorem}[1][]%
{\begin{trivlist}\item[]\refstepcounter{theoremcntr}%
{\bfseries Theorem~\thetheoremcntr%
  \ifthenelse{\equal{#1}{}}{}{ (#1)}.
}}%
{\hfill$\Box$\end{trivlist}}

\newcounter{propositioncntr}
{\begin{trivlist}\item[]\refstepcounter{propositioncntr}%
{\bfseries Proposition~\thepropositioncntr%
  \ifthenelse{\equal{#1}{}}{}{ (#1)}.
}}%
{\hfill$\Box$\end{trivlist}}

\newcounter{lemmacntr}
{\begin{trivlist}\item[]\refstepcounter{lemmacntr}%
{\bfseries Lemma~\thelemmacntr%
  \ifthenelse{\equal{#1}{}}{}{ (#1)}.
}}%
{\hfill$\Box$\end{trivlist}}

\newcommand{\restrict}[2]{\left.#1\right|_{#2}}

\newcommand{\imk}{{1 - \kappa(s) z}}

\newcommand{\va}{{\sigma_a^2}}
\newcommand{\vb}{{\sigma_b^2}}
\newcommand{\hva}{{{\hat \sigma_a}^2}}
\newcommand{\hvb}{{{\hat \sigma_b}^2}}

\newcommand{\mth}{{m_\theta}}

\newcommand{\vth}{{\sigma_\theta^2}}

\newcommand{\calN}{\mathcal{N}}

\newcommand{\cent}[1]{\makebox(0,0){#1}}
\newcommand{\pos}[2]{\makebox(0,0)[#1]{#2}}

\newcommand{\knownBox}{\cent{\rule{1.75\unitlength}{1.75\unitlength}}}

\newfloat{algbox}{t}{alg}
\floatname{algbox}{Algorithm}

\definecolor{gray}{rgb}{0.5, 0.5, 0.5}

\newcommand{\gray}{\color{gray}}

\begin{document}
\title{Model-Predictive Control with NUP Priors}
\author{Raphael~Keusch and Hans-Andrea Loeliger
\thanks{R.~Keusch was, and H.-A.~Loeliger is, with the Department of Information
Technology and Electrical Engineering, ETH Zurich, 8092 Zurich, Switzerland
(e-mail: raphael.keusch@bluewin.ch, loeliger@isi.ee.ethz.ch).}%
}%

\markboth{XXXXXX}%
{Keusch \MakeLowercase{\textit{et al.}}: Model-Predictive Control with New NUV Priors}

\maketitle

\begin{abstract}
  Normals with unknown variance (NUV) 
  and, more generally, normals with unknown parameters (NUP)
  can represent many useful priors 
  including $\text{L}_p$ norms and other sparsifying priors, 
  and they blend well with linear-Gaussian models 
  and Gaussian message passing algorithms.
  In this paper, we elaborate on recently proposed NUP representations 
  of half-space constraints, box constraints, and finite-level constraints.
  We then demonstrate the use of such NUP representations 
  for exemplary applications in model predictive control 
  with a variety of constraints on the input, the output, or the internal state
  of the controlled system. In such applications, the computations boil down 
  to iterations of Kalman-type forward-backward recursions,
  with a complexity (per iteration) that is linear in the planning horizon.
  In consequence, this approach can handle long planning horizons, 
  which distinguishes it from the prior art.
  For nonconvex constraints, this approach has no claim to optimality, 
  but it is empirically very effective.
\end{abstract}


\section{Introduction}

\IEEEPARstart{N}{ormal} priors with unknown variance (NUV priors) 
are a central idea of sparse Bayesian learning~\cite{Tipping2001, Tipping2003, Wipf2004, Wipf2008}
and closely related to variational representations of cost functions
and iteratively reweighted least-squares methods~\cite{Bach2012, Daubechies2010, Loeliger2018}.
The point of such priors is the computational compatibility with linear Gaussian models. 
The primary use of such priors has been to encourage sparsity,
in applications including sparse input estimation~\cite{Loeliger2016, Zalmai2017}, 
localized event detection~\cite{zalmai_blind_2016, Zalmai2017a, wadehn_model-based_2020}, 
outlier removal~\cite{wadehn_outlier-insensitive_2016, Wadehn2019}, 
sparse least squares~\cite{Loeliger2018}, control~\cite{Bruderer2015, Hoffmann2017a}, 
and imaging~\cite{Ma2020, ma2022}.

A next step was made by the binarizing NUP prior 
(normal with unknown parameters)
recently proposed 
in~\cite{keusch2021binaryNUV, Marti21multiuser},
which may be viewed as consisting of two internal NUV priors.
Such composite NUV priors offer many additional possibilities, some of which 
are proposed and explored in this paper.

Specifically, in this paper, we propose\footnote{
The first write-ups of these new composite NUV priors 
are~\cite{keusch2021binaryNUVext,keusch2021boxPreprint}
which have not otherwise been published; see also \cite{Keusch2022PhD}.
}
and explore
a NUP prior
to enforce half-plane constraints,
and we generalize the mentioned binarizing NUP prior to $M$-level priors with \mbox{$M>2$}.

We then demonstrate the use of such NUP priors 
for a variety of control problems
with constraints on inputs, outputs, and states.
In such applications, the computations amount to iterating Kalman-type
forward-backward recursions, with simple closed-form updates 
of the NUP parameters in between.
The computational complexity of each such iteration is linear in ``time''
(i.e., in the planning horizon); in consequence, this approach 
can handle long planning horizons (with high temporal resolution), 
which distinguishes it from the prior art.

Using NUP priors to express constraints comes with 
the qualification
that constraints can be enforced only for variables that are 
deterministic functions of the controlling input.
Therefore, in this paper, we restrict ourselves 
to cases where the controlled system 
is entirely deterministic with known initial state.

The related literature of constrained optimization is vast. 
Numerous methods have been developed in the broader field of constrained convex optimization with 
linear inequality constraints---most notably the
projected Newton method~\cite{bertsekas1982projected,kim2010tackling},
the projected gradient method~\cite{rosen1960gradient}, 
the interior-point method~\cite{wright1997primal}, 
and the active set method~\cite{stark1995bounded}.
Generally speaking, the computational complexity of
these methods scales much faster than linearly 
with the number of constraints.
Methods such as \cite{MuehleDAndrea2019}
solve this problem in some cases, but not in general.

Discrete-level constraints generically results in NP-hard problems.
Finding the optimal solution to such problems using exhaustive enumeration is thus 
limited to short planning horizons~\cite{FSQS09}.
Another naive approach is to first solve the unconstrained problem 
and then project the solution to the feasible set;
unfortunately, most often, the obtained solution is far from optimal.
Tree-search algorithms with branch-and-bound methods such as 
sphere decoding~\cite{hassibi2005sphere, KaGK16} may help, but their complexity is 
still exponential in the planning horizon. By contrast, the approach of this 
paper offers excellent empirical performance with linear complexity in the 
planning horizon.

The paper is structured as follows.
The idea of statistical models with NUP priors 
is briefly reviewed in Section~\ref{sec:SysNUV}.
The proposed approach to control is developed in Section~\ref{sec:ProposedApproach}.
The NUP representation for half-space constraints 
and discrete-level constraints 
are derived and discussed in Sections \ref{sec:HalfSpace}
and \ref{sec:DiscreteLevel}, respectively.

Section~\ref{sec:app} demonstrates the application of the proposed method 
to a variety of exemplary constrained control problems including 
bounded-error control, 
binary and ternary control,
and minimal-time race track control.
In a companion paper~\cite{keusch2023ContrSysTechn}, the approach of this paper is 
applied to 
a real-world control problem in power electronics.

The following notation is used.
The Gaussian probability density function in $x$ with mean $m$ and 
covariance matrix $V$ is denoted by $\Normal{x; m, V}$. 
Equality of functions up to a scale factor is denoted by $\propto$.



\section{Models with NUP Priors: A Brief Review}
\label{sec:SysNUV}

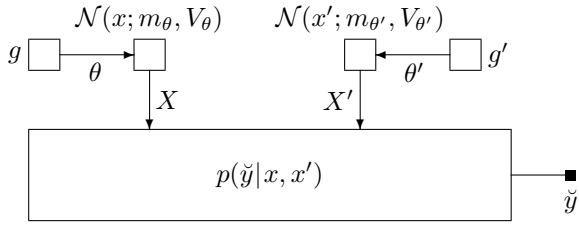
\begin{figure}
\centering
\setlength{\unitlength}{0.8mm}
\begin{picture}(91,35.5)(0,10)
%
\put(17.5,35){\framebox(5,5){}}  \put(20,41.25){\pos{cb}{$\calN(x; m_\theta, V_\theta)$}}
 \put(0,35){\framebox(5,5){}}    \put(-1.25,37.5){\pos{cr}{$g$}}
 \put(5,37.5){\vector(1,0){12.5}}  \put(11,36.25){\pos{ct}{$\theta$}}
\put(20,35){\vector(0,-1){10}}   \put(21,30){\pos{cl}{$X$}}
\put(52.5,35){\framebox(5,5){}}  \put(55,41.25){\pos{cb}{$\calN(x'; m_{\theta'}, V_{\theta'})$}}
 \put(70,35){\framebox(5,5){}}   \put(76.25,37.5){\pos{cl}{$g'$}}
 \put(70,37.5){\vector(-1,0){12.5}}  \put(64,36.25){\pos{ct}{$\theta'$}}
\put(55,35){\vector(0,-1){10}}   \put(54,30){\pos{cr}{$X'$}}
\put(0,10){\framebox(80,15){$p(\breve y \cond x, x')$}}
\put(80,17.5){\line(1,0){10}}  \put(90,15){\pos{ct}{$\breve y$}}
\put(90,17.5){\knownBox}
\end{picture}

\vspace{2mm}

\caption{\label{fig:NUP:Splitting}%
Factor graph of system model (\ref{eqn:NUV:Sys}) 
with NUP priors (\ref{eqn:NUP}) 
and fixed observation(s) \mbox{$\breve Y = \breve y$}.}
\end{figure}

Consider a statistical model
with latent variables $X$ and $X'$, observable(s)%
\footnote{The accent of $\breve Y$ is for compatibility with Section~\ref{sec:ProposedApproach}.}
$\breve Y$,
and (improperly scaled)
joint probability density function
\begin{equation} \label{eqn:NUV:Sys}
p(\breve y, x, x'; \theta, \theta') = 
p(\breve y \cond x, x') p(x; \theta) p(x'; \theta')
\end{equation}
with unknown parameters $\theta$ and $\theta'$,
where
\begin{equation} \label{eqn:NUP}
p(x; \theta) = \calN(x; m_\theta, V_\theta) g(\theta),
\end{equation}
i.e., $p(x; \theta)$ is Gaussian in $x$ (up to a scale factor) 
with mean $m_\theta$, 
variance (or covariance matrix) $V_\theta$, 
and nonnegative scale factor $g(\theta)$
all depending on $\theta$,
and likewise for $p(x'; \theta') g'(\theta')$,
as illustrated in Fig.~\ref{fig:NUP:Splitting}.
The factor $g(\theta)$ in (\ref{eqn:NUP}) 
may be viewed as a (improper) prior on $\theta$,
and (\ref{eqn:NUP}) may thus be viewed as a (improper) joint prior on 
$X$ and $\theta$.

We further assume
that the likelihood 
$p(\breve y \cond x, x')$ is jointly Gaussian in $x$ and $x'$ (up to a scale factor). 
It follows that, for fixed $\breve y$ and fixed $\theta$ and $\theta'$, 
the posterior $p(x,x' \cond \breve y; \theta, \theta')$ is jointly Gaussian in $x$ and $x'$.

The point of the NUP representation (\ref{eqn:NUP})
is to reduce 
estimation in the model (\ref{eqn:NUV:Sys}) to iterations 
of estimation in a Gaussian model.
The effect of the NUP prior (\ref{eqn:NUP}) depends 
on how the unknown parameters $\theta$ and $\theta'$ are estimated,
as will be detailed below.

The generalization of (\ref{eqn:NUV:Sys}) to $n$ variables $X_1,\dots,X_n$
will be obvious throughout.

\subsection{Joint MAP Estimation with Alternating Maximization}
\label{sec:SysNUV:AM}

In this approach, we aim to compute (for fixed $\breve Y = \breve y$)
\begin{IEEEeqnarray}{rCl}
\mleft( \hat x, \hat x' \mright) 
& = &
\argmax_{x,\, x'} \max_{\theta,\, \theta'} p(\breve y, x, x'; \theta, \theta') 
\label{eqn:NUV:Sys:JointMAP}\\
& = &  \argmax_{x,\, x'} p(\breve y \cond x, x') \rho(x) \rho'(x')
\end{IEEEeqnarray}
where 
\begin{equation} \label{eqn:NUV:Sys:effprior}
\rho(x) \eqdef \max_{\theta} p(x; \theta)
\end{equation}
and $\rho'(x')$ (defined analogously)
are the effective (possibly improper%
\footnote{
In (\ref{eqn:NUV:Sys:effprior}), we use $\rho$ instead of $p$ 
to emphasize that it may not be normalizable.}%
) priors on $X$ and $X'$, respectively.
An obvious approach to compute the estimate (\ref{eqn:NUV:Sys:JointMAP})
is to repeat the following two steps 
for $i=1,2,3,\ldots$ until convergence
(beginning with some initial $\theta^{(0)}$ and $(\theta')^{(0)}$):
\begin{enumerate}
\item
For fixed $\theta=\theta^{(i-1)}$ and $\theta'=\theta^{(i-1)}$, compute 
\begin{IEEEeqnarray}{rCl}
\IEEEeqnarraymulticol{3}{l}{
\big( x^{(i)}, (x')^{(i)} \big)
 =  \argmax_{x,\, x'} p(\breve y, x, x'; \theta, \theta')
}\\
& = & \argmax_{x,\, x'} p(\breve y \cond x, x') \calN(x; m_\theta, V_\theta) \calN(x'; m_{\theta'}, V_{\theta'}),
      \IEEEeqnarraynumspace
      \label{eqn:NUV:Sys:GaussianMAP}
\end{IEEEeqnarray}
which is a Gaussian MAP/MMSE estimate.
\item
For fixed $x=x^{(i)}$ and $x'=(x')^{(i)}$, compute
\begin{equation}
\Big( \theta^{(i)}, (\theta')^{(i)} \Big) 
= \argmax_{\theta,\, \theta'} p(\breve y, x, x'; \theta, \theta'),
\end{equation}
which splits into
\begin{equation}  \label{eqn:NUV:Sys:JointMAP:updateNUP}
\theta^{(i)} = \argmax_{\theta} p(x; \theta)
\end{equation}
and likewise for $\theta'$.
\end{enumerate}

\subsection{Type-II MAP Estimation%
\protect\footnote{in the sense of \cite{Tipping2001, Wipf2004}} 
with Expectation Maximization}
\label{sec:SysNUV:EM}

In this approach, we aim to first compute (for fixed $\breve Y = \breve y$)
\begin{equation} \label{eqn:NUV:Sus:TypeIIMAP}
\big( \hat \theta, \hat \theta' \big) 
=
\argmax_{\theta,\, \theta'} \int_{x}\int_{x'} p(\breve y, x, x'; \theta, \theta')\, dx'\, dx, 
\end{equation}
after which $\hat x$ and $\hat x'$ are computed by (\ref{eqn:NUV:Sys:GaussianMAP}).
The maximization in (\ref{eqn:NUV:Sus:TypeIIMAP}) 
is carried out by expectation maximization (EM) with hidden variables $X$ and $X'$,
i.e., by iterating
\begin{equation}  \label{eqn:NUV:Sys:JointEM}
\big( \theta^{(i)}, (\theta')^{(i)} \big)
= \argmax_{\theta,\, \theta'} \E\mleft[ \log p(\breve y, X, X'; \theta, \theta')\rule{0em}{2ex} \mright],
\end{equation}
where the expectation is computed 
with respect to 
\begin{IEEEeqnarray}{rCl}
\IEEEeqnarraymulticol{3}{l}{
p\big( x, x' \cond \breve y; \theta^{(i-1)}, (\theta')^{(i-1)} \big)
}\nonumber\\\quad
 & \propto &
p(\breve y \cond x, x') p(x; \theta^{(i-1)}) p(x'; (\theta')^{(i-1)}),
\IEEEeqnarraynumspace
\end{IEEEeqnarray}
i.e., in a jointly Gaussian setting as in (\ref{eqn:NUV:Sys:GaussianMAP}).
The maximization (\ref{eqn:NUV:Sys:JointEM}) splits into
\begin{equation} \label{eqn:NUV:Sys:UpdateEM}
\theta^{(i)} = \argmax_{\theta} \E\mleft[ \log p(X; \theta)\rule{0em}{2ex} \mright]
\end{equation}
and likewise for $\theta'$.
The computation of (\ref{eqn:NUV:Sys:UpdateEM}) 
boils down to the computation of the mean and the variance (or the covariance matrix) of $X$
for fixed $\theta=\theta^{(i-1)}$ and $\theta'=(\theta')^{(i-1)}$,
cf.\ Section~\ref{sec:disc:TypeII}.

In this approach, there is no counterpart to the effective prior (\ref{eqn:NUV:Sys:effprior}).

\subsection{Mixed MAP Estimation}
\label{sec:SysNUV:Mixed}

We will sometimes mix the methods of Sections \ref{sec:SysNUV:AM} and~\ref{sec:SysNUV:EM}
by updating some parameters using (\ref{eqn:NUV:Sys:UpdateEM}) 
while updating some other parameters using (\ref{eqn:NUV:Sys:JointMAP:updateNUP}). 
We thus effectively maximize some compromise between 
a cost function as in (\ref{eqn:NUV:Sys:JointMAP}) and a cost function as in (\ref{eqn:NUV:Sus:TypeIIMAP}).

\section{Proposed Approach}
\label{sec:ProposedApproach}

\subsection{System Model and Examples}
\label{sec:SysApproachExamples}

Recall the standard linear state space model
\begin{IEEEeqnarray}{rCl}  \IEEEyesnumber \phantomsection \label{eqn:lssm:detLSSM}  \IEEEyessubnumber* 
  x_k &=& A x_{k-1} + B u_k \\
  y_k &=& C x_k
\end{IEEEeqnarray}
with time index $k\in\{ 1, 2, 3, \ldots\}$,
input $u_k \in \Reals^L$, 
state $x_k \in \Reals^N$,
output $y_k \in \Reals^H$, 
and matrices $A \in \Reals^{N \times N}$,
$B \in \Reals^{N \times L}$,
and $C \in \Reals^{H \times N}$.
We assume that the input $u_1, u_2, \ldots$ 
can be used to control the system.

For a given initial state $x_0$, 
a given planning horizon%
\footnote{%
The extension to model-predictive control with a receding planning horizon
is straightforward, cf.\ Fig.~\ref{fig:versus_optimal_both_short_horizon}
and~\cite{keusch2023ContrSysTechn}.
} 
$K$,
a given target 
$\breve y = \bma \breve y_1, \ldots, \breve y_K \ema$,
and optionally a given target state $\breve x_K$,
we wish to compute a control $u = \bma u_1, \ldots, u_K \ema$
that minimizes some cost function subject to constraints
on $u$, $y = \bma y_1, \ldots, y_K \ema$, and $x = \bma x_0, \ldots, x_K \ema$, 
as illustrated by the following examples.
In these first examples, for ease of exposition, the input and the output are scalar 
and there are no costs or constraints on the states;
more examples (without these restrictions) are discussed in 
Section~\ref{sec:app} and in \cite{keusch2023ContrSysTechn}.

\begin{example} \label{example:prob:classicControl}
  Classical linear-quadratic control problem:
  \begin{IEEEeqnarray}{rCl}
    \hat u = \argmin_u \|y(u) - \breve y\|^2 + \alpha \|u\|^2
  \end{IEEEeqnarray}
  for some given $\alpha \in\Reals, \alpha > 0$.
\end{example}

\begin{example} Squared fitting error with sparse input (the LASSO problem~\cite{Tibshirani1996}): 
  \begin{IEEEeqnarray}{rCl}
    \hat u = \argmin_u \|y(u) - \breve y\|^2 + \alpha \|u\|_1
  \end{IEEEeqnarray}
  for some given $\alpha \in\Reals, \alpha > 0$.
\end{example}

\begin{example} \label{example:prob:binaryControl}
  Squared fitting error with binary input: 
  \begin{IEEEeqnarray}{rCl} \IEEEyesnumber \phantomsection \label{eqn:prob:binaryControlOptProb}  \IEEEyessubnumber*
    \hat u &=& \argmin_u \|y(u) - \breve y\|^2 \quad \text{s.t.} \label{eqn:prob:binaryControlOptProbQuadPart}\\
         u_k &\in& \{0, 1\} \quad \text{or} \quad u_k \in \{-1, +1\}, \quad k \in \{1, \dots, K\}. \label{eqn:prob:binaryControlOptProbBinConstr}
  \end{IEEEeqnarray}
\end{example}

\begin{example} 
$\text{L}_1$ fitting error and bounded input: 
\begin{IEEEeqnarray}{rCl} \IEEEyesnumber \IEEEyessubnumber*
  \hat u &=& \argmin_u \|y(u) - \breve y\|_1 \quad \text{s.t.} \\
       u_k &\in& [a, b], \quad k \in \{1, \dots, K\}
\end{IEEEeqnarray}
for given $a, b \in \Reals$.
\end{example}

\begin{example} Bounded fitting error and sparse input level switches: 
  \begin{IEEEeqnarray}{rCl} \IEEEyesnumber \IEEEyessubnumber*
    \hat u &=& \argmin_u \|\Delta u\|_0 \quad \text{s.t.} \\
         && | y_k(u) - \breve y_k | \leq b, \quad k \in \{1, \dots, K\}
  \end{IEEEeqnarray}
  for some given $b \in \Reals, b > 0$, and $\Delta u \eqdef \bma u_2 \!-\! u_1, \ldots, u_K \!-\! u_{K-1} \ema$.
\end{example}
\noindent
Note that Example~\ref{example:prob:classicControl} is a classical control problem,
which is well-known to be solvable by Kalman-type forward-backward recursions,
with complexity linear in $K$ (cf. Section~\ref{sec:IAKE}).

The essence of this paper is that all these problems 
(and many other combinations of constraints and cost functions
on inputs, outputs, and states)
can be efficiently solved%
---exactly in the convex case, otherwise approximately---%
by an iterative algorithm, where each iteration solves a statistical 
estimation problem that is essentially equivalent to (some variation of)
Example~\ref{example:prob:classicControl}.

\subsection{The Statistical Model}

The equivalence of Example~\ref{example:prob:classicControl}
with MAP/MMSE estimation in a linear-Gaussian model
is standard:
\begin{IEEEeqnarray}{rCl}
\hat u
& = & \argmin_u \|y(u) - \breve y\|^2 + \alpha \|u\|^2 \\
& = & 
      \argmax_{u} \exp\mleft( - \frac{\| y(u) - \breve y \|^2}{2\sigma^2} \mright)
      \exp\mleft( - \frac{\| u \|^2}{2\sigma^2\alpha^{-1}} \mright) 
      \IEEEeqnarraynumspace\label{eqn:StatModel:JustGaussian}\\
& = & \argmax_{u} p( \breve y \cond u) p(u)
\end{IEEEeqnarray}
with arbitrary $\sigma^2>0$, 
where
\begin{equation}
p(u) \eqdef \prod_{k=1}^K p(u_k)
\end{equation}
with
\begin{equation}
p(u_k) \eqdef \Normal{u_k; 0, \sigma^2/\alpha}
\end{equation}
and
\begin{equation}
p( \breve y \cond u) \eqdef \prod_{k=1}^K p\big( \breve y_k \cond y_k(u) \big)
\end{equation}
with
\begin{IEEEeqnarray}{rCl}
p\big( \breve y_k \cond y_k(u) \big) 
& \eqdef & \Normal{\breve y_k; y_k(u), \sigma^2} \IEEEeqnarraynumspace\\
& = &  \Normal{y_k(u); \breve y_k, \sigma^2}.
              \label{eqn:GaussModel:LikelihoodFactor}
\end{IEEEeqnarray}

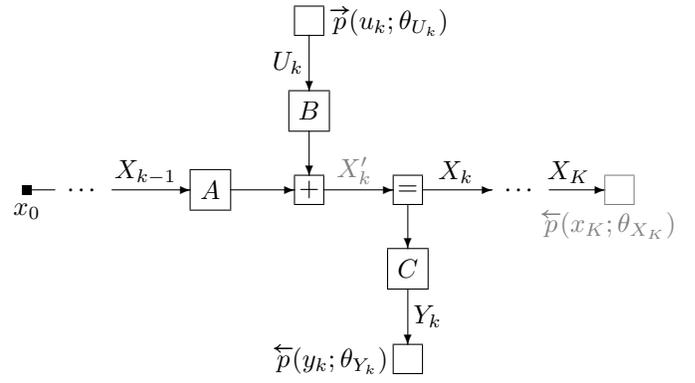
\begin{figure}
\centering
\setlength{\unitlength}{0.75mm}
\begin{picture}(117,65)(7.5,0)

\put(57.5,60){\framebox(5,5){}}   \put(63.5,62.5){\pos{cl}{$\msgf{p}{}(u_k; \theta_{U_k})$}}
\put(60,60){\vector(0,-1){10}}   \put(59,55){\pos{cr}{$U_k$}}
\put(56.5,43){\framebox(7,7){$B$}}
\put(60,43){\vector(0,-1){8}}

\put(10,32.5){\knownBox}    \put(10,30){\pos{ct}{$x_0$}}
\put(10,32.5){\line(1,0){5}}    
\put(20,32.5){\cent{\ldots}}
\put(25,32.5){\vector(1,0){14}}    \put(31,33.5){\pos{cb}{$X_{k-1}$}}
\put(39,29){\framebox(7,7){$A$}}  
\put(46,32.5){\vector(1,0){11.5}}
\put(57.5,30){\framebox(5,5){$+$}}
\put(62.5,32.5){\vector(1,0){12.5}}   \put(68,33.5){\pos{cb}{\gray$X_k'$}}
\put(75,30){\framebox(5,5){$=$}}
\put(80,32.5){\vector(1,0){12.5}}   \put(86,33.75){\pos{cb}{$X_k$}}
\put(97.5,32.5){\cent{\ldots}}
\put(102.5,32.5){\vector(1,0){10}}  \put(106,33.75){\pos{cb}{$X_K$}}
{\gray
\put(112.5,30){\framebox(5,5){}}   \put(113,29){\pos{ct}{$\msgb{p}{}(x_K; \theta_{X_K})$}}
}

\put(77.5,30){\vector(0,-1){8}}
\put(74,15){\framebox(7,7){$C$}}
\put(77.5,15){\vector(0,-1){10}}  \put(78.5,10){\pos{cl}{$Y_k$}}
\put(75,0){\framebox(5,5){}}    \put(74,2.5){\pos{cr}{$\msgb{p}{}(y_k; \theta_{Y_k})$}}
\end{picture}
\vspace{-1mm}
\caption{\label{fig:lssm:jointDensity}%
Factor graph of the model~(\ref{eqn:GenStatModel})
and (\ref{eqn:GenStatModelwithFinal}).}
\end{figure}

We now generalize this linear Gaussian model to
\begin{equation} \label{eqn:GenStatModel}
  p(y, x, u; \theta) \!\propto\!  \prod_{k=1}^K \msgf{p}{}(u_k;\theta_{U_k})
  \msgb{p}{}(y_k;\theta_{Y_k})  \bigg\rvert_{\text{(\ref{eqn:lssm:detLSSM})}}  
\end{equation}
or to
\begin{IEEEeqnarray}{rCl} \label{eqn:GenStatModelwithFinal}
  p(y, x, u; \theta) \!\propto\! \mleft( \prod_{k=1}^K \msgf{p}{}(u_k;\theta_{U_k})
  \msgb{p}{}(y_k;\theta_{Y_k}) \! \! \mright) \msgb{p}{}(x_K; \theta_{X_K}) \bigg\rvert_{\text{(\ref{eqn:lssm:detLSSM})}}  
  \nonumber\\
\end{IEEEeqnarray}
(as illustrated in Fig.~\ref{fig:lssm:jointDensity}) 
with 
\begin{equation}
\theta \eqdef (\theta_{U_1}, \dots, \theta_{U_K}, \theta_{Y_1} \dots , \theta_{Y_K}, \theta_{X_K}),
\end{equation}
where 
$\msgf{p}{}(u_k;\theta_{U_k})$, $\msgb{p}{}(y_k;\theta_{Y_k})$, 
and $\msgb{p}{}(x_K; \theta_{X_k})$
comprise NUP representations as in (\ref{eqn:NUP}):
for fixed $\theta_{U_k} \eqdef (\msgf{m}{U_k}, \msgf{V}{U_k})$,
\begin{equation} \label{eqn:p:msgfUk}
\msgf{p}{}(u_k;\theta_{U_k}) \propto \Normal{u_k; \msgf{m}{U_k}, \msgf{V}{U_k}},
\end{equation}
for fixed $\theta_{Y_k} \eqdef (\msgb{m}{Y_k}, \msgb{V}{Y_k})$,
\begin{equation} \label{eqn:p:msgbYk}
\msgb{p}{}(y_k;\theta_{Y_k}) \propto \Normal{y_k; \msgb{m}{Y_k}, \msgb{V}{Y_k}},
\end{equation}
and likewise for the optional factor $\msgb{p}{}(x_K; \theta_{X_K})$.
The arrows in (\ref{eqn:p:msgfUk}) and (\ref{eqn:p:msgbYk}) 
refer to the direction of the corresponding edges in Fig.~\ref{fig:lssm:jointDensity}
and distinguish between priors and likelihoods.
Note that $\msgb{m}{Y_k}$
subsumes the target $\breve y_k$ in (\ref{eqn:GaussModel:LikelihoodFactor}).

In the statistical model (\ref{eqn:GenStatModel}) or (\ref{eqn:GenStatModelwithFinal}),
inputs, outputs, and states are random variables and therefore
denoted by capital letters $U_k$, $Y_k$, and $X_k$, respectively.

\subsection{Multiple Inputs, Outputs, and State Constraints}
\label{sec:MultInOut}

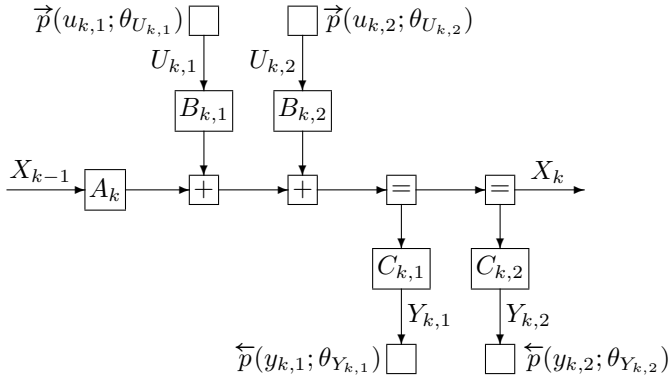
\begin{figure}
\centering
\setlength{\unitlength}{0.75mm}
\begin{picture}(112,65)(7.5,0)

\put(40,60){\framebox(5,5){}}     \put(39,62.5){\pos{cr}{$\msgf{p}{}(u_{k,1}; \theta_{U_{k,1}})$}}
\put(42.5,60){\vector(0,-1){10}}   \put(41.5,55){\pos{cr}{$U_{k,1}$}}
\put(37.5,43){\framebox(10,7){$B_{k,1}$}}
\put(42.5,43){\vector(0,-1){8}}

\put(57.5,60){\framebox(5,5){}}   \put(63.5,62.5){\pos{cl}{$\msgf{p}{}(u_{k,2}; \theta_{U_{k,2}})$}}
\put(60,60){\vector(0,-1){10}}   \put(59,55){\pos{cr}{$U_{k,2}$}}
\put(55,43){\framebox(10,7){$B_{k,2}$}}
\put(60,43){\vector(0,-1){8}}

\put(7.5,32.5){\vector(1,0){14}}   \put(13.5,33.5){\pos{cb}{$X_{k-1}$}}
\put(21.5,29){\framebox(7,7){$A_k$}}  
\put(28.5,32.5){\vector(1,0){11.5}}
\put(40,30){\framebox(5,5){$+$}}
\put(45,32.5){\vector(1,0){12.5}}
\put(57.5,30){\framebox(5,5){$+$}}
\put(62.5,32.5){\vector(1,0){12.5}}   
\put(75,30){\framebox(5,5){$=$}}
\put(80,32.5){\vector(1,0){12.5}}   
\put(92.5,30){\framebox(5,5){$=$}}
\put(97.5,32.5){\vector(1,0){12.5}}   \put(103.5,33.75){\pos{cb}{$X_k$}}

\put(77.5,30){\vector(0,-1){8}}
\put(72.5,15){\framebox(10,7){$C_{k,1}$}}
\put(77.5,15){\vector(0,-1){10}}  \put(78.5,10){\pos{cl}{$Y_{k,1}$}}
\put(75,0){\framebox(5,5){}}    \put(74,2.5){\pos{cr}{$\msgb{p}{}(y_{k,1}; \theta_{Y_{k,1}})$}}

\put(95,30){\vector(0,-1){8}}
\put(90,15){\framebox(10,7){$C_{k,2}$}}
\put(95,15){\vector(0,-1){10}}  \put(96,10){\pos{cl}{$Y_{k,2}$}}
\put(92.5,0){\framebox(5,5){}}   \put(99,2.5){\pos{cl}{$\msgb{p}{}(y_{k,2}; \theta_{Y_{k,2}})$}}

\end{picture}
\vspace{3mm}
\caption{\label{fig:lssm:SplitInOut}%
Multiple inputs and outputs as in Section~\ref{sec:MultInOut}.}
\end{figure}

The following variations of the system model (\ref{eqn:lssm:detLSSM}) 
are used in some of the examples 
in Section~\ref{sec:app} and in \cite{keusch2023ContrSysTechn}.

First, the generalization 
to time-varying matrices $A_k$, $B_k$, and $C_k$ (e.g., from linearizing a nonlinear model) 
is immediate.

Second, multidimensional inputs $U_k$ and outputs $Y_k$ can sometimes be split 
into lower-dimensional (preferably scalar) inputs $U_{k,1}, U_{k,2}, \ldots$
or outputs $Y_{k,1}, Y_{k,2}, \ldots$
(with individual NUP factors $\msgf{p}{}(u_{k,1}; \theta_{U_{k,1}})$ etc.)
as illustrated in Fig.~\ref{fig:lssm:SplitInOut}.

Third, additional (virtual) outputs $Y_{k,\ell}$ 
with pertinent vectors or matrices $C_{k,\ell}$ and 
NUP factors $\msgb{p}{}(y_{k,\ell}; \theta_{k,\ell})$
can be used to impose constraints on linear functions 
of~$X_k$. 

\subsection{Iterative Augmented Kalman Estimation (IAKE)}
\label{sec:IAKE}

For fixed NUP parameters $\theta$, 
all inputs $U_k$, states $X_k$, and outputs $Y_k$ in the statistical model
(\ref{eqn:GenStatModel}) or (\ref{eqn:GenStatModelwithFinal})
are jointly Gaussian,
and the MAP estimate of any subset of these variables 
coincides with their (posterior) mean. 
For the joint estimation of all these variables and $\theta$, 
we will use the following iterative algorithm,
which implements both methods of Section~\ref{sec:SysNUV}.

Starting from an initial guess $\hat \theta^{(0)}$, the algorithm repeats the 
following two steps for $i=1,2,3,\ldots,$ until convergence (or for a sufficiently large number of 
iterations):

\begin{algbox}
  \centering
  \noindent

\framebox[0.49\textwidth]{
  \begin{minipage}{0.45\textwidth}
    \newcounter{saveequationcntr}
    \setcounter{saveequationcntr}{\value{equation}}
    \setcounter{equation}{0}
    \renewcommand{\theequation}{F.\arabic{equation}}
    
\noindent
\vskip0.2em
\rule{0pt}{3ex}The algorithm consists of a forward recursion followed by a backward recursion.
The former is a standard Kalman filter, but the latter is not quite standard.

  \vspace{0.2cm}
  \textbf{Forward recursion:} \\
  Initialize with $\msgf{m}{X_0}=x_0$ and $\msgf{V}{X_0}=0^{N\times N}$. Then,
  for $k = 1, 2, \dots, K$, compute  
  \vspace{-0.2cm}
  \begin{IEEEeqnarray}{rCl}
    \msgf{m}{X_k'} & = &  A \msgf{m}{X_{k-1}} + B \msgf{m}{U_k}  \\
    \msgf{V}{X_k'} &=& A \msgf{V}{X_{k-1}} A^{\T} + B \msgf{V}{U_k} B^{\T} \IEEEeqnarraynumspace \\
    \msgf{m}{X_k} & = &  \msgf{m}{X_k'} + \msgf{V}{X_k'} E_k \\
    \msgf{V}{X_k} &=& F_{k} \msgf{V}{X_k'}  \\
    %
    \noalign{\noindent where \vspace{2\jot}}
    E_{k} &=& C^{\T} G_{k} (\msgb{m}{Y_{k}} -  C \msgf{m}{X_{k}'}) \\
    F_{k} &=& I_N - \msgf{V}{X_{k}'} C^{\T} G_{k}  C \\
    G_{k} &=& (\msgb{V}{Y_{k}} + C \msgf{V}{X_{k}'} C^{\T})^{-1}.
  \end{IEEEeqnarray}

  \setcounter{equation}{0}
  \renewcommand{\theequation}{B.\arabic{equation}}
  
  \textbf{Backward recursion:} \\
  If $\msgb{p}{}(x_K; \theta_{X_K}) = 1$, 
  initialize with $\tilde W_{X_K}\!= 0_{N \times N}$ and $\tilde \xi_{X_K}\!= 0_N$; 
  otherwise
  \vspace{-0.2cm}
  \begin{IEEEeqnarray}{rCl}
    \tilde W_{X_K} &=& ( \msgf{V}{X_K} + \msgb{V}{X_K} )^{-1}  \label{eqn:MBF:FinalStateVarInv}\\ 
    \tilde \xi_{X_K} &=& \tilde W_{X_K} (\msgf{m}{X_K} - \msgb{m}{X_K}). \hspace{1.3em} \IEEEeqnarraynumspace
  \end{IEEEeqnarray}
  For $k = K, K-1, \dots, 1$, compute 
  \vspace{-0.1cm}
  \begin{IEEEeqnarray}{rCl}
    \tilde{\xi}_{X_k'} &=& F_k^{\T} \tilde{\xi}_{X_k} - E_k  \\
    \tilde{W}_{X_k'} &=& F_k^{\T} \tilde{W}_{X_k} F_k + C^{\T} G_k C  \hspace{2.em} \IEEEeqnarraynumspace \\
    \tilde{\xi}_{X_{k-1}} &=& A^{\T} \tilde{\xi}_{X_k'}  \\
    \tilde{W}_{X_{k-1}} &=& A^{\T} \tilde{W}_{X_k'} A.
  \end{IEEEeqnarray}

  \setcounter{equation}{0}
  \renewcommand{\theequation}{P.\arabic{equation}}
  
  \textbf{Posterior quantities (= estimates):} \\
  The posterior means and variances for $k \in \{1, \dots, K\}$ are given by
  \vspace{-0.2cm}
  \begin{IEEEeqnarray}{rCl}
    m_{U_k} &=& \msgf{m}{U_k} - \msgf{V}{U_k} B^{\T} \tilde \xi_{X_k'} \label{eqn:algo:mU}\\
    V_{U_k} &=& \msgf{V}{U_k} - \msgf{V}{U_k} B^{\T} \tilde{W}_{X_k'} B \msgf{V}{U_k} \label{eqn:algo:VU} 
                 \IEEEeqnarraynumspace\\
    m_{X_k} &=& \left ( \msgf{m}{X_k} - \msgf{V}{X_k} \tilde \xi_{X_k} \right ) \label{eqn:algo:mX}\\
    V_{X_k} &=& \left ( \msgf{V}{X_k} - \msgf{V}{X_k} \tilde{W}_{X_k} \msgf{V}{X_k} \right )  \label{eqn:algo:VX}\\
    m_{Y_k} &=& C m_{X_k} \label{eqn:algo:mY}\\
    V_{Y_k} &=& C V_{X_k} C^{\T}. \label{eqn:algo:VY}
  \end{IEEEeqnarray}
  \vspace{-0.3cm}
  \setcounter{equation}{\value{saveequationcntr}}
\end{minipage}
}

\caption{\label{alg:MBFMP}%
Step~\ref{item:step_1} of IAKE implemented by 
MBF message passing with input estimation assembled from \cite{Loeliger2016}.
In many applications, only a subset of~(\ref{eqn:algo:mU})--(\ref{eqn:algo:VY}) needs to be computed.
}
\end{algbox}

\begin{enumerate}
\item \label{item:step_1}
  For fixed $\theta = \hat \theta^{(i-1)}$, compute for $k \in \{1, \dots, K\}$
  \begin{enumerate}
    \item the posterior means $m_{U_k}^{(i)}$ (and, if necessary, the posterior variances 
              $V_{U_k}^{(i)}$) of $U_k$, 
    \item the posterior means $m_{Y_k}^{(i)}$ (and, if necessary, the posterior variances 
              $V_{Y_k}^{(i)}$) of $Y_k$. 
  \end{enumerate}
\item \label{item:step_2}
  From these means and variances, determine new parameters $\theta^{(i)}$ using 
  implementations of (\ref{eqn:NUV:Sys:JointMAP:updateNUP}) and/or (\ref{eqn:NUV:Sys:UpdateEM})
  as in Tables~\ref{tab:prob:UpdateRulesNUV} and~\ref{tab:psum:UpdateRulesCNUV}
  (which will be discussed below).
\end{enumerate}
The output of the algorithm 
(i.e., the desired control sequence) 
is the final estimate 
$\hat u_1 = m_{U_1}$, \ldots, $\hat u_K = m_{U_K}$.

\begin{table*}
  \setcounter{saveequationcntr}{\value{equation}}
  \setcounter{equation}{0}
  \renewcommand{\theequation}{T\ref{tab:prob:UpdateRulesNUV}.\arabic{equation}}
  \newcommand{\putEqnNo}[1]{%
     \begin{minipage}{3em}
     \vspace{-5ex}
     \begin{equation}\label{#1}
     \end{equation}%
     \end{minipage}%
  }

  \normalsize
  \centering
  \renewcommand{\arraystretch}{1.5}
  \begin{tabular}[t]{@{}llclr@{}}
  \toprule
   Prior & Use Case && \multicolumn{1}{l}{Update Rules} & \\
  \toprule
  
  \makecell[l]{$\text{L}_1$ (Laplace) \cite{Bach2012, Loeliger2018}} & 
  sparsity &&
  $\msgf{V}{X} = \gamma^{-1} |m_X|$ & \putEqnNo{eqn:TableNUV:L1}
  \\
  \midrule

  \makecell[l]{$\text{L}_p$ \cite{Bach2012, ma2022}} & 
  various &&
  $\displaystyle \msgf{V}{X} = \frac{|m_X|^{2-p}}{\gamma p}$ & \putEqnNo{eqn:TableNUV:Lp}
  \\
  \midrule

  \makecell[l]{Smoothed $\text{L}_1$/Huber \cite{Loeliger2018, ma2022}} & 
  \makecell[l]{outlier-insensitive fitting} &&
  $\displaystyle \msgf{V}{X} = \max \left \{ r^2, \frac{|m_X|}{\gamma} \right \}$
  & \putEqnNo{eqn:TableNUV:Huber}
  \\
  \midrule

  \makecell[l]{Plain NUV \cite{Loeliger2016}} & 
  \makecell[l]{sparsity} &&
  $\msgf{V}{X} = V_X + m_X^2$
  & \putEqnNo{eqn:TableNUV:plainNUV}
  \\
  \midrule

  \makecell[l]{Smoothed plain NUV \cite{Ma2020, ma2022}} & 
  \makecell[l]{outlier-insensitive fitting} &&
  \makecell[l]{ 
  $\msgf{V}{X} = \max \left \{ r^2, m_X^2 \right \} \quad $ or\\
  $\msgf{V}{X} = \max \left \{ r^2, V_X + m_X^2 \right \}$
  }
  & 
  \makecell[l]{
  \putEqnNo{eqn:TableNUV:SmoothedPlainNUVAM} \\
  \rule{0em}{3ex}\putEqnNo{eqn:TableNUV:SmoothedPlainNUVEM}
  }
  \\
  \bottomrule
  \\[-3ex]
  \end{tabular}
\caption{\label{tab:prob:UpdateRulesNUV}%
Update rules for 
some
basic NUV priors 
(with parameters $\gamma$ and $r^2$), 
cf.\ the cited references. 
The mean $\msgf{m}{X}$ remains zero.
}
\setcounter{equation}{\value{saveequationcntr}}
\end{table*}

\begin{table*}
  \setcounter{saveequationcntr}{\value{equation}}
  \setcounter{equation}{0}
  \renewcommand{\theequation}{T\ref{tab:psum:UpdateRulesCNUV}.\arabic{equation}}
  \newcommand{\putEqnNo}[1]{%
     \begin{minipage}{3em}
     \vspace{-5ex}
     \begin{equation}\label{#1}
     \end{equation}%
     \end{minipage}%
  }

  \normalsize
  \centering
  \renewcommand{\arraystretch}{1.5}
  \begin{tabular}[t]{@{}llclr@{}}
  \toprule
   Prior & Constraint && \multicolumn{1}{l}{Update Rules} & \\
  \toprule
  \multirow{3}{*}{\makecell[l]{\vspace{0mm}\\ Hinge loss\\(Section~\ref{sec:halfSpacePrior})}} & 
  $x \geq a$ && 
  \makecell[l]{ 
    $\displaystyle \msgf{V}{X} = \frac{\left | m_X - a \right|}{\gamma} $ \vspace{0.1cm} \\ 
    $\msgf{m}{X} = a + \left | m_X - a \right|$ \vspace{0.15cm}
  } 
  & \putEqnNo{eqn:TableNUP:geq}
  \\
  \cline{2-5}
    & $x \leq a$ && 
    \rule{-3pt}{6.5ex} 
  \makecell[l]{
    $\displaystyle \msgf{V}{X} = \frac{\left | m_X - a \right|}{\gamma} $ \vspace{0.1cm} \\ 
    $\msgf{m}{X} = a - \left | m_X - a \right|$  \vspace{0.1cm}
    } 
    & \putEqnNo{eqn:TableNUP:leq}
    \\
    \midrule
  \makecell[l]{Vapnik loss\\(Section~\ref{sec:boxPrior})} & 
  $a \leq x \leq b$ && 
  \makecell[l]{ 
    $\displaystyle \msgf{V}{X} =  \frac{1}{\gamma}\left ( \frac{1}{\left | m_X - a \right|} + \frac{1}{\left | m_X - b \right|} \right )^{-1}$  \vspace{0.2cm} \\
    $\displaystyle \msgf{m}{X} = \gamma\msgf{V}{X}  \left ( \frac{a}{\left | m_X - a \right|} + \frac{b}{\left | m_X - b \right|} \right)$ \vspace{0.1cm} \\
  } 
  & \putEqnNo{eqn:TableNUP:box}
  \\
  \midrule

  \makecell[l]{Binarizing prior\\(\cite{keusch2021binaryNUV} and Section~\ref{sec:DiscreteLevel})} & 
  $x \in \{a, b\}$ && 
  \makecell[l]{
    $\displaystyle \msgf{V}{X}= \left ( \! \frac{1}{ V_X \!+\! (m_{X} \!-\! a )^2 } \!+\! \frac{1}{V_X \!+\! (m_{X} \!-\! b )^2 } \! \right )^{\!-1}$  \vspace{0.1cm}\\ 
    $\displaystyle \msgf{m}{X} = \msgf{V}{X} \! \left ( \! \frac{a}{ V_X \!+\! (m_{X} \!-\! a )^2 } \!+\! \frac{b}{V_X \!+\! (m_{X} \!-\! b )^2 } \! \right )$   \vspace{0.1cm}
  }
  & \putEqnNo{eqn:TableNUP:Bin}
  \\
  \bottomrule
  \\[-3ex]
  \end{tabular}
\caption{\label{tab:psum:UpdateRulesCNUV}%
Update rules for 
the NUP priors
of Sections \ref{sec:HalfSpace} and~\ref{sec:DiscreteLevel}.
}
\setcounter{equation}{\value{saveequationcntr}}
\end{table*}

Note that Step~\ref{item:step_1} operates with a standard linear Gaussian model.
In consequence, the required means and variances can be computed 
by Kalman-type recursions or, equivalently, by forward-backward Gaussian 
message passing, with a complexity that is linear in $K$.

A preferred such algorithm is the
Modified Bryson--Frazier (MBF) smoother~\cite{Bierman1977}
augmented with input signal estimation as in~\cite{Bruderer2015, Loeliger2016}.
This algorithm does not require to invert any $N\times N$ matrices 
(except perhaps (\ref{eqn:MBF:FinalStateVarInv}))
and it is numerically stable. 
For the convenience 
of the reader, 
the algorithm is concisely stated 
as Algorithm~\ref{alg:MBFMP}.

The modular approach of \cite{Loeliger2016} makes it easy 
to adapt this algorithm (and other Kalman-type algorithms)
to variations in the setting (e.g., as in Section~\ref{sec:MultInOut}).

To guarantee constraint satisfaction, 
the algorithm of this section may need to be repeated
a few times as discussed in Section~\ref{sec:ConstraintsChecking}.

\subsection{Tabulated Update Rules for Step~2 of IAKE}

Tables~\ref{tab:prob:UpdateRulesNUV} and~\ref{tab:psum:UpdateRulesCNUV} 
give the update rules 
according to (\ref{eqn:NUV:Sys:JointMAP:updateNUP}) or (\ref{eqn:NUV:Sys:UpdateEM})
for the NUP parameters of some generic scalar variable $X$
(which can be applied both to $X=U_k$ 
(with $\msgf{m}{X} = \msgf{m}{U_k}$ and $\msgf{V}{X} = \msgf{V}{U_k}$)
and to $X=Y_k$
(with $\msgf{m}{X} = \msgb{m}{Y_k}$ and $\msgf{V}{X} = \msgb{V}{Y_k}$).

The update rules in Table~\ref{tab:prob:UpdateRulesNUV} are not new 
and stated here only for the sake of completeness.
For vector versions of Table~\ref{tab:prob:UpdateRulesNUV}, we refer to 
\cite{LgMLSP2023}.

Table~\ref{tab:psum:UpdateRulesCNUV} will be discussed 
in Sections \ref{sec:HalfSpace} and~\ref{sec:DiscreteLevel}.

\subsection{Outer Loop for Constraint Satisfaction}
\label{sec:ConstraintsChecking}

For the priors of Table~\ref{tab:psum:UpdateRulesCNUV}
to properly act as constraints, the scale parameter $\gamma$
in (\ref{eqn:TableNUP:geq})--(\ref{eqn:TableNUP:box}) must be sufficiently large
(as will be discussed in Section~\ref{sec:HalfSpace}),
or the variance $\sigma^2$ of Gaussian factors as in (\ref{eqn:StatModel:JustGaussian})
must be sufficiently large
(as will be discussed in Section~\ref{sec:DiscreteLevel}).
However, suitable minimal values of these parameters are not known a priori,
and increasing $\gamma$ slows down IAKE 
(by requiring more iterations).

This issue can be addressed as follows.
Beginning with initial values (e.g.,~1) for these parameters,
repeat the following steps until all constraints are satisfied:
{\renewcommand{\theenumi}{\alph{enumi}}
\begin{enumerate}
\item Run IAKE (Section~\ref{sec:IAKE}) to convergence.
\item Check if all constraints are satisfied. If yes, we are done.
\item Increase the pertinent scale factors $\gamma$ 
      or/and the variance $\sigma^2$ of Gaussian factors (e.g., by a factor of~2).
\end{enumerate}
Some specific examples will be discussed in Section~\ref{sec:app}.
}

In practice (i.e., excluding adversarial problem statements)
the required number of such outer iterations
seems to be limited to some small number such as~10.

\subsection{Why Deterministic Systems?}

NUP priors per se belong to statistical models. 
However, handling constraints by NUP priors as in
Table~\ref{tab:psum:UpdateRulesCNUV} 
constrains the estimate, not the actual values, 
of the corresponding variables.

In the setting of this paper, 
constraints on the actual values can thus be enforced
only for variables that are (deterministic) functions
of the control sequence $u$ (including $u$ itself). 
Therefore, in this paper, we restrict ourselves 
to deterministic systems as in Section~\ref{sec:SysApproachExamples},
which allows us to impose constraints 
on all variables (inputs, outputs, and states).


\section{NUP Priors for Half-Space Constraints and Box Constraints}
\label{sec:HalfSpace}

\begin{figure}
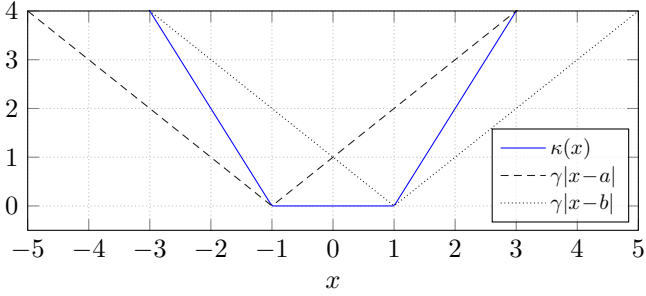
   
  \centering
    \includestandalone{inc_pdf/boxConstraintCostFunc}
    \caption{%
    Cost function~(\ref{eqn:box:CostFunctionBoxPrior}) 
    for $a=-1, b=1$, and $\gamma = 1$.
    }
    \label{fig:box:CostBoxPrior} 
\end{figure}

This section is about 
(\ref{eqn:TableNUP:geq})--(\ref{eqn:TableNUP:box})
in Table~\ref{tab:psum:UpdateRulesCNUV}.
Let
\begin{equation} \label{eqn:box:CostFunctionBoxPrior}
\kappa(x) \eqdef \gamma \big( |x-a| + |x-b| - |b-a| \big)
\end{equation}
be the Vapnik loss function (cf.\ Fig.~\ref{fig:box:CostBoxPrior})
and let
\begin{equation} \label{eqn:box:priorModelEffective}
p_\text{V}(x) \propto \exp\!\big( {-}\kappa(x) \big)
\end{equation}
be the associated (normalizable) prior.
The idea is to use (\ref{eqn:box:priorModelEffective})
(with sufficiently large $\gamma$) to enforce 
an estimate $\hat x$ with $a \leq \hat x \leq b$.
In a second step, we obtain a half-space constraint
by a suitable limit $a\rightarrow -\infty$ or $b\rightarrow\infty$.

In this section, we use joint MAP estimation as in 
Section~\ref{sec:SysNUV:AM} 
with NUP representations as in (\ref{eqn:NUV:Sys:effprior}).

\subsection{NUP Representation of Vapnik Loss}
\label{sec:boxPrior}

The Laplace prior has the (well-known) NUV representation
\begin{IEEEeqnarray}{rCl}  \label{eqn:box:laplaceNUVrep}
  \exp \! \big (\!-\! \gamma |x| \big ) = \max_{\sigma^2} \Normal{x; 0, \sigma^2} \tilde g(\sigma^2),
\end{IEEEeqnarray}
with $\gamma > 0$ and 
$\tilde g(\sigma^2) \eqdef \sqrt{2 \pi \sigma^2} e^{-\gamma^2 \sigma^2/2}$~\cite{Loeliger2018},
and the maximizing variance in~(\ref{eqn:box:laplaceNUVrep}) is
\begin{IEEEeqnarray}{rCl} \label{eqn:box:laplaceSigmaHat}
  \hat \sigma^2 = \argmax_{\sigma^2} \Normal{x;0, \sigma^2} \tilde g(\sigma^2) = \frac{|x|}{\gamma}. 
\end{IEEEeqnarray}
Thus (\ref{eqn:box:priorModelEffective}) can be written as
\begin{IEEEeqnarray}{rCl}
  p_\text{V}(x) & \propto & \exp \! \big (\!-\! \gamma |x-a| \big ) \cdot \exp \! \big (\!-\! \gamma |x-b| \big ) \\
  &=& \max_{\va} \Normal{x;a, \va} \tilde g(\va) \cdot \max_{\vb} \Normal{x;b, \vb} \tilde g(\vb).
  \IEEEeqnarraynumspace \label{eqn:box:NUVrep}
\end{IEEEeqnarray}
Using (\ref{eqn:ProductGaussians}), we then obtain the NUP representation
\begin{equation} \label{eqn:VapnikNUP}
p_\text{V}(x) \propto \max_{\theta} \Normal{x; m_\theta, \sigma^2_\theta} g(\theta)
\end{equation}
with $\theta \eqdef (\va, \vb)$,
\begin{IEEEeqnarray}{rCl} \label{eqn:box:boxPriorMeanVar}
  \vth = \left ( \frac{1}{\va} + \frac{1}{\vb} \right)^{-1} \text{~and~~~} \mth = \vth \left ( \frac{a}{\va} + \frac{b}{\vb} \right),
   \IEEEeqnarraynumspace
\end{IEEEeqnarray}
and
\begin{equation}
g(\theta) = \Normal{a-b; 0, \sigma_a^2 + \sigma_b^2} \tilde g(\va) \tilde g(\vb).
\end{equation}

From~(\ref{eqn:box:laplaceSigmaHat}), the maximizing variances 
in (\ref{eqn:box:NUVrep}) and~(\ref{eqn:VapnikNUP})
are
\begin{IEEEeqnarray}{rCl}
  \hva = \frac{|x-a|}{\gamma} \quad \text{and} \quad \hvb = \frac{|x-b|}{\gamma}. \label{eqn:box:boxOptimalVars}
\end{IEEEeqnarray}
Plugging 
(\ref{eqn:box:boxOptimalVars}) 
into~(\ref{eqn:box:boxPriorMeanVar}) yields
\begin{IEEEeqnarray}{rCl}
  \vth &=& \left ( \frac{\gamma}{\left | x - a \right|} + \frac{\gamma}{\left | x - b \right|} \right )^{-1} \label{eqn:box:VXfBox}\\
  \mth &=& \gamma\vth  \left ( \frac{a}{\left | x - a \right|} + \frac{b}{\left | x - b \right|} \right), \label{eqn:box:mXfBox}
\end{IEEEeqnarray}
which is (\ref{eqn:TableNUP:box})
(in slightly different notation).

\subsection{Box Constraint: Single-Variable Analysis} 
\label{sec:box:trivialExample}

\begin{figure}
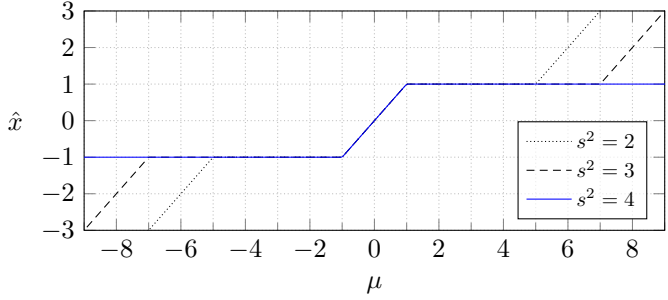

  \centering
  \includestandalone{inc_pdf/boxConstraintChar}
  \caption{\label{fig:box:boxConstChar}%
  Estimate~(\ref{eqn:box:Analysis:EstX}) for $a=-1$, $b=1$, $\gamma=1$,
  and different values of $s^2$. 
 }
\end{figure}

We next study the effect of (\ref{eqn:VapnikNUP}) as a box constraint.
Consider a statistical model with latent variable $X$, 
observation $\breve Y = \breve y$, 
and joint probability density function
\begin{equation} \label{eqn:box:Analysis:Model}
p(\breve y, x; \theta) = p(\breve y \cond x) \Normal{x;\mth, \vth} g(\theta). 
\end{equation}
We further assume
\begin{IEEEeqnarray}{rCl} \label{eqn:box:GaussianLikelihood}
p(\breve y \cond x) \propto \Normal{x;\mu, s^2}.
\end{IEEEeqnarray}
(Note that this is the setting of Section~\ref{sec:SysNUV}
without~$X'$.)
Clearly, joint MAP estimation of $\theta$ and $X$ yields
\begin{IEEEeqnarray}{rCl}
\hat x & = & \argmax_{x} \max_\theta p(\breve y, x; \theta) \\
& = &  \argmax_{x} p(\breve y \cond x) p_\text{V}(x).
     \label{eqn:box:Analysis:EstX}
\end{IEEEeqnarray}
The estimate (\ref{eqn:box:Analysis:EstX}) as a function of $\mu$
is plotted in Fig.~\ref{fig:box:boxConstChar}.
We observe that, for given $\mu$ and $\gamma > 0$ and 
with sufficiently large $s^2$, the estimate (\ref{eqn:box:Analysis:EstX}) is 
indeed restricted to $[a,b]$.
Quantitatively, we have
\begin{theorem} \label{theorem:box:hardConstraint}
Let $a<b$. The estimate~(\ref{eqn:box:Analysis:EstX}) satisfies $a \leq \hat x \leq b$ 
if and only if
\begin{IEEEeqnarray}{rCl} \label{eqn:box:HardCondition}
  s^2 > \begin{cases}
  0  & \text{if $a \leq \mu \leq b$}, \\
  \min \left \{ \frac{|a-\mu|}{2\gamma}, 
                \frac{|b - \mu|}{2\gamma} \right \} &
  \text{otherwise}.
  \end{cases} 
\vspace{-2.5ex}
\end{IEEEeqnarray}
\end{theorem}
The proof is not difficult and given in Appendix~\ref{sec:apds:ProofBox}.
It is thus obvious that 
the constraint $a \leq \hat x \leq b$
can be enforced by choosing $\gamma$ to be sufficiently large.

Finally, it is not hard to see that 
alternatingly maximizing (\ref{eqn:box:Analysis:Model})
over $\theta$ and over $x$ (as in Section~\ref{sec:SysNUV:AM})
will 
converge to (\ref{eqn:box:Analysis:EstX})
except for very unlucky initializations such as  $\va = 0$ or $\vb = 0$.

\subsection{NUP Representation of Hinge Loss}
\label{sec:halfSpacePrior}

\begin{figure}
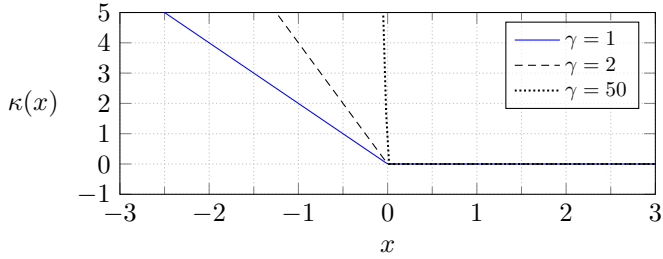

  \centering
  \includestandalone{inc_pdf/CostBoxPriorPlanePlus}
  \caption{\label{fig:hsc:costFuncPlus}%
  Cost function~(\ref{eqn:hsc:costFuncPlus}) 
  for $a = 0$ and different values of $\gamma$.
  }
  \end{figure}

Taking the limit $b\rightarrow\infty$ in~(\ref{eqn:box:CostFunctionBoxPrior}) 
yields 
\begin{IEEEeqnarray}{rCl} \label{eqn:hsc:costFuncPlus}
  \kappa(x) &=& 
  \begin{cases} 2 \gamma (a - x) & \text{if $x < a$},  \\ 
              0 &\text{otherwise},
  \end{cases} \IEEEeqnarraynumspace 
\end{IEEEeqnarray}
which is illustrated in Fig.~\ref{fig:hsc:costFuncPlus}.
Taking the limit $b \rightarrow \infty$ in~(\ref{eqn:box:VXfBox}) 
and~(\ref{eqn:box:mXfBox}) yields
\begin{IEEEeqnarray}{rCl}
  \sigma_\theta^2 
  & = &  \frac{|x - a|}{\gamma} \label{eqn:hsc:mXfHscRight}\\
  m_\theta
  & = & a + |x - a|, 
  \label{eqn:hsc:VXfHscRight} \IEEEeqnarraynumspace
\end{IEEEeqnarray}
which is (\ref{eqn:TableNUP:geq})
(in slightly different notation).

Likewise, taking the limit $b\rightarrow -\infty$ in~(\ref{eqn:box:CostFunctionBoxPrior}) 
yields 
\begin{IEEEeqnarray}{rCl} \label{eqn:hsc:costFuncMinus}
  \kappa(x) &=& 
  \begin{cases} 0 & \text{if $x < a$}, \\ 
    2 \gamma (x - a) &\text{otherwise},
\end{cases} \IEEEeqnarraynumspace
\end{IEEEeqnarray}
and taking the limit 
$b \rightarrow -\infty$ in~(\ref{eqn:box:VXfBox}) and~(\ref{eqn:box:mXfBox})
yields
\begin{IEEEeqnarray}{rCl}
   \sigma_\theta^2 &=& 
  \frac{|x - a|}{\gamma} \label{eqn:hsc:mXfHscLeft}\\
  m_\theta &=& 
  a - |x - a|.
   \label{eqn:hsc:VXfHscLeft}
\end{IEEEeqnarray}
which is (\ref{eqn:TableNUP:leq}).

\subsection{Half-Space Constraint: Single-Variable Analysis}
\label{sec:hsc:trivialExample}

\begin{figure}
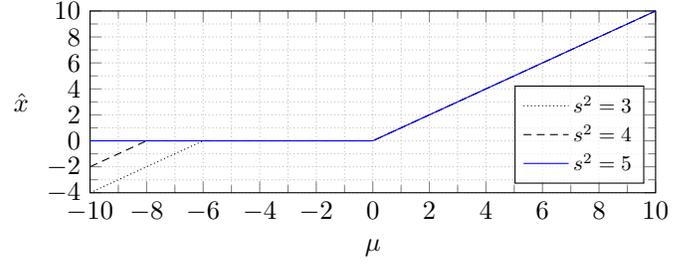

  \centering
  \includestandalone{inc_pdf/halfPlaneConstraintChar}
  \caption{\label{fig:hsc:halfPlaneConstraintChar}%
  Estimate~(\ref{eqn:hsc:jointMAPSimplified}) for $a=0$ 
  and different values of $s^2$. 
  }
\end{figure}

Let
\begin{equation} \label{eqn:hsc:priorPlus}
\rho_+(x) \eqdef {\exp}\big( {-}\kappa(x) \big)
\end{equation}
with $\kappa(x)$ as in (\ref{eqn:hsc:costFuncPlus}).
Analyzing the effect of (\ref{eqn:hsc:priorPlus}) as a half-space constraint
amounts to a simplified version of Section~\ref{sec:box:trivialExample}.
The estimate
\begin{equation} \label{eqn:hsc:jointMAPSimplified}
\hat x = \argmax_x  p(\breve y \cond x) \rho_+(x)
\end{equation}
with $p(\breve y \cond x)$ as in (\ref{eqn:box:GaussianLikelihood})
is illustrated in Fig.~\ref{fig:hsc:halfPlaneConstraintChar}.
We observe that, for any fixed $\mu$ and $\gamma$ and 
sufficiently large $s^2$, 
the estimate~(\ref{eqn:hsc:jointMAPSimplified}) 
indeed satisfies $\hat x \geq a$. 
Quantitatively, we have
\begin{theorem} \label{theorem:hsc:hardConstraint}
  The estimate~(\ref{eqn:hsc:jointMAPSimplified}) satisfies
  $\hat x \geq a$ if and only if 
  \begin{IEEEeqnarray}{rCl} \label{eqn:hsc:HardConditionHPRight}
    s^2 >  \begin{cases}
      0 & \text{if $\mu \geq a$}, \\
      \frac{|a-\mu|}{2\gamma} &
      \mathrm{otherwise}.
      \end{cases} 
  \vspace{-2ex}
  \end{IEEEeqnarray}
\end{theorem}
(The proof of Theorem~\ref{theorem:hsc:hardConstraint} is easily obtained as a 
suitably simplified version of the proof of Theorem~\ref{theorem:box:hardConstraint}.)
It is thus obvious that the constraint $\hat x \geq a$
can be enforced by choosing $\gamma$ to be sufficiently large.

\section{NUP Priors for Discrete-Level Constraints}
\label{sec:DiscreteLevel}

This section is about (\ref{eqn:TableNUP:Bin})
in Table~\ref{tab:psum:UpdateRulesCNUV}: 
we discuss and generalize 
the composite-NUP prior for enforcing $x\in \{ a, b\}$
that was proposed in~\cite{keusch2021binaryNUV}.
%
%
%
This prior is given by
\begin{IEEEeqnarray}{rCl} \label{eqn:disc:TwoLvlPrior}
p(x; \theta) \eqdef \Normal{x; a, \va} \Normal{x; b, \vb}
\end{IEEEeqnarray}
with $\theta \eqdef (\va, \vb)$.
It turns out that this prior strongly prefers $X$ to lie in $\{ a, b\}$.
The detailed working of this binarizing effect depends
on how the unknown variances $\theta$ are determined,
as will be discussed below.

Using (\ref{eqn:ProductGaussians}), (\ref{eqn:disc:TwoLvlPrior}) can be written as 
\begin{IEEEeqnarray}{rCl} \label{eqn:disc:PriorWithHyperPrior}
p(x; \theta) =  \Normal{x; \mth, \vth} g(\theta)
\end{IEEEeqnarray}
with $\mth$ and $\vth$ as in~(\ref{eqn:box:boxPriorMeanVar}), 
and
\begin{IEEEeqnarray}{rCl} \label{eqn:disc:HyperPrior}
g(\theta) =  \Normal{a-b; 0, \va + \vb}.
\end{IEEEeqnarray}

\subsection{Joint MAP Estimation: Effective Prior}

\begin{figure}
  \centering
  \includestandalone[]{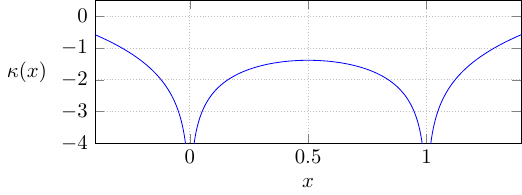}
  \caption{The cost function~(\ref{eqn:disc:CostFunctionBinPrior}) for $a=0$ and $b=1$.}
  \label{fig:disc:costFuncBinPrior}%
  \end{figure}

We first consider estimating $\theta$ as in Section~\ref{sec:SysNUV:AM}
(but estimating $\theta$ as in Section~\ref{sec:SysNUV:EM} works much better,
as will be discussed below).
It is easily seen that
\begin{equation} \label{eqn:binary:maximizingVars}
\argmax_{\theta} p(x; \theta) 
= (\hva, \hvb) 
= \left( (x - a)^2,  (x - b)^2 \right).
\end{equation}
Plugging this into (\ref{eqn:disc:PriorWithHyperPrior}) 
yields the effective prior (\ref{eqn:NUV:Sys:effprior}) 
\begin{IEEEeqnarray}{rCl}
\rho(x) & = & \max_{\theta} p(x; \theta) \label{eqn:disc:binPriorEffectiveGen}\\
& \propto & \frac{1}{|x-a| \cdot |x-b|}
             \label{eqn:disc:binPriorEffective}
\end{IEEEeqnarray}
The associated cost function 
\begin{IEEEeqnarray}{rCl}
  \kappa(x) = -\log p(x) =  \log|x-a| +  \log|x-b| + \operatorname{const.} 
  \IEEEeqnarraynumspace
  \label{eqn:disc:CostFunctionBinPrior} 
\end{IEEEeqnarray}
is illustrated in Fig.~\ref{fig:disc:costFuncBinPrior}.
It is obvious that this prior 
strongly favors $X$ to lie in $\{a, b\}$.
  
Plugging (\ref{eqn:binary:maximizingVars}) into
$\Normal{x; \mth, \vth}$
with (\ref{eqn:box:boxPriorMeanVar}) yields
\begin{IEEEeqnarray}{rCl}
  \sigma_\theta^2 &=& \left ( \frac{1}{  (x - a )^2 } + \frac{1}{ (x - b )^2 } \right )^{-1} \label{eqn:disc:VXfAM}\\ 
  m_\theta &=& \sigma_\theta^2  \left ( \frac{a}{  (x - a )^2 } + \frac{b}{ (x - b )^2 } \right ), \label{eqn:disc:mXfAM}
\end{IEEEeqnarray}
which would be the update rules as in Table~\ref{tab:psum:UpdateRulesCNUV}
(but are not actually stated there).

\subsection{Joint MAP Estimation: Single-Variable Analysis}
\label{sec:disc:JointMAP:Analysis}

Analogously to Section~\ref{sec:box:trivialExample},
we next study the effect of (\ref{eqn:disc:binPriorEffectiveGen}) as a binarizing constraint.
Consider a statistical model with latent variable $X$, 
observation $\breve Y = \breve y$, 
and joint probability density function
\begin{equation} \label{eqn:binary:Analysis:Model}
p(\breve y, x; \theta) = p(\breve y \cond x) \Normal{x;\mth, \vth} g(\theta). 
\end{equation}
We further assume
\begin{IEEEeqnarray}{rCl} \label{eqn:binary:GaussianLikelihood}
p(\breve y \cond x) \propto \Normal{x;\mu, s^2}.
\end{IEEEeqnarray}
Clearly, joint MAP estimation of $\theta$ and $X$ yields
\begin{IEEEeqnarray}{rCl}
\hat x & = & \argmax_{x} \max_\theta p(\breve y, x; \theta) \\
& = &  \argmax_{x} p(\breve y \cond x) \rho(x).
     \label{eqn:binary:Analysis:EstX}
\end{IEEEeqnarray}

\begin{figure}
  \centering
  \includestandalone[]{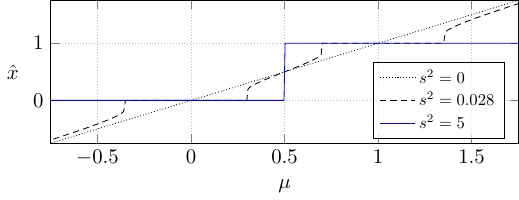}
  \caption{\label{fig:hatXJointMAP}%
  The estimate (\ref{eqn:binary:Analysis:EstX})
  for $a=0$ and $b=1$. 
  }
\end{figure}

The estimate (\ref{eqn:binary:Analysis:EstX}) as a function of $\mu$
is plotted in Fig.~\ref{fig:hatXJointMAP}.
We observe that for given $\mu$ and a sufficiently large $s^2$, 
the estimate discretizes, i.e., $\hat x \in \{a, b\}$.
Quantitatively, we have
\begin{theorem}\label{theorem:AM:scalarLocalMaxima}
The function
\begin{IEEEeqnarray}{rCl}
x \mapsto \frac{\Normal{x; \mu, s^2}}{|x-a|\cdot|x-b|}
\end{IEEEeqnarray}
has no local maximum
(other than the global maxima at $x=a$ and $x=b$) 
if and only if
\begin{IEEEeqnarray}{rCl} \label{eqn:disc:AM:CondNoLocal}
s^2 > s^2_{\mathrm{AM}},
\end{IEEEeqnarray}
where $s^2_{\mathrm{AM}}$ depends on $\mu, a$ and $b$.
\end{theorem}

\begin{figure}
\centering
\includestandalone[]{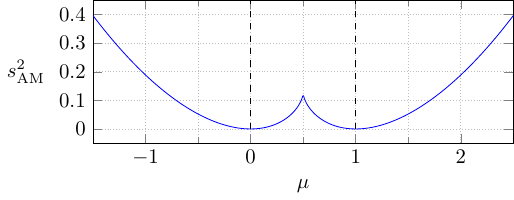}
\caption{\label{fig:BoundAM}%
The value of $s^2_{\text{AM}}$ in (\ref{eqn:disc:AM:CondNoLocal})
as a function of $\mu$, 
for $a=0$ and $b=1$.}
\end{figure}

The proof of Theorem~\ref{theorem:AM:scalarLocalMaxima} 
(including the definition of $s^2_\text{AM}$) is lengthy 
and omitted here
but can be found in \cite{keusch2021binaryNUVext} and \cite{Keusch2022PhD}. 
Since $s^2_{\text{AM}}$ is the only real root of a cubic polynomial,
a closed-form expression for $s^2_{\text{AM}}$ exists, 
but it is cumbersome. 
However, $s^2_{\text{AM}}$ is easily computed numerically.
The value of $s^2_{\text{AM}}$ as a function of $\mu$
is plotted in Fig.~\ref{fig:BoundAM}. 
For example, $s_{\text{AM}}^2 = 0.028$ for $\mu = 0.3$, $a=0$, and $b=1$
(cf. Fig.~\ref{fig:hatXJointMAP}).

If (\ref{eqn:disc:AM:CondNoLocal}) holds, 
alternatingly maximizing (\ref{eqn:binary:Analysis:Model})
over $\theta$ and over $x$ (as in Section~\ref{sec:SysNUV:AM})
will converge to $\hat x = a$ or to $\hat x = b$ 
(except if $x$ is unluckily initialized to the unavoidable
local minimum between $a$ and~$b$). 
However, there is a catch:
depending on the initialization, 
this alternating maximization may get 
trapped into the wrong maximum.
This pitfall is avoided by estimation 
as in Section~\ref{sec:disc:TypeII}.

\subsection{Type-II MAP Estimation: Update Rule}
\label{sec:disc:TypeII}

We next consider estimating $\theta$ as in Section~\ref{sec:SysNUV:EM},
which turns out to work much better.
We first work out the update rule (\ref{eqn:NUV:Sys:UpdateEM}):
\begin{IEEEeqnarray}{rCl} 
\theta^{(i)}
& = & \argmax_{\theta} \E\mleft[ \log p(X; \theta)\rule{0em}{2ex} \mright] \\
& = & \argmax_{\va,\, \vb} \E\mleft[ \log\mleft( \Normal{X; a, \va} \Normal{X; b, \vb}\rule{0em}{2.5ex}\mright) \mright], 
      \IEEEeqnarraynumspace
\end{IEEEeqnarray}
which splits into
\begin{IEEEeqnarray}{rCl}
\mleft( \va \mright)^{(i)} 
& = & \argmax_{\va} \E\mleft[ \log\Normal{X; a, \va} \mright] \\
& = & \argmin_{\va} \mleft(
      \frac{1}{2}\log(\va)
      + \frac{1}{2\va}\EE{ \left( X - a \right)^2 }
      \mright). 
      \IEEEeqnarraynumspace\label{eqn:apds:EMderStep1} 
\end{IEEEeqnarray}
and likewise for $\vb$.
Setting the derivative 
with respect to $\va$ 
to zero yields
\begin{IEEEeqnarray}{rCl}
\mleft( \va \mright)^{(i)} 
& = &  \EE{ \left( X - a \right)^2 }  \label{eqn:apds:EMderStep2} \\
&=& \EE{X^2} - \EE{X}^2 + \EE{X}^2 - 2a \EE{X} + a^2 \IEEEeqnarraynumspace \\ 
&=& \Var{X} + (\EE{X} - a)^2
\end{IEEEeqnarray}
and likewise for $\vb$.
Plugging these updates for $\va$ and $\va$ into
$\Normal{x; \mth, \vth}$ 
with (\ref{eqn:box:boxPriorMeanVar}) yields
\begin{IEEEeqnarray}{rCl}
  \sigma_\theta^2 &=&  \left ( \frac{1}{ V_X + (m_X - a )^2 } + \frac{1}{V_X + (m_X - b )^2 } \right )^{-1} 
                \IEEEeqnarraynumspace\label{eqn:disc:VXfEM}\\ 
  m_\theta &=&  \sigma_\theta^2 \left ( \frac{a}{ V_X + (m_X - a )^2 } + \frac{b}{V_X + (m_X - b )^2 } \right ) 
                \IEEEeqnarraynumspace\label{eqn:disc:mXfEM}
\end{IEEEeqnarray}
with $V_X \eqdef \Var{X}$ and $m_X \eqdef \EE{X}$,
which is (\ref{eqn:TableNUP:Bin}).

\subsection{Type-II MAP Estimation: Single-Variable Analysis}
\label{sec:disc:TypeII:Analysis}

Consider the statistical model 
(\ref{eqn:binary:Analysis:Model}) and (\ref{eqn:binary:GaussianLikelihood}) 
as in Section~\ref{sec:disc:JointMAP:Analysis}.
Using (\ref{eqn:disc:VXfEM}) and (\ref{eqn:disc:mXfEM}), 
we now estimate $X$ by expectation maximization as in 
Section~\ref{sec:SysNUV:EM}. 

\begin{figure}
  \centering
  \includestandalone[]{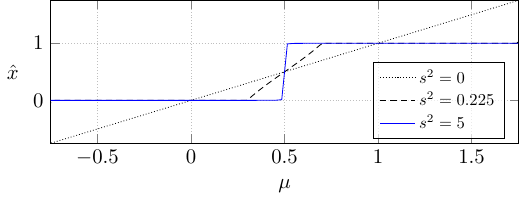}
  \caption{\label{fig:hatXTypeII}%
  The estimate of Section~\ref{sec:disc:TypeII:Analysis}
  for $a=0$ and $b=1$. 
  }  
\end{figure}

\begin{figure}
\centering
\includestandalone[]{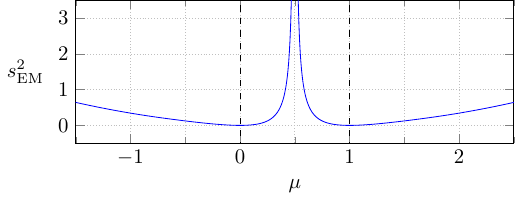}
\caption{\label{fig:BoundEM}%
The value of $s_{\text{EM}}^2$ in (\ref{eqn:disc:scalarEMtheorem:Condsa}) and (\ref{eqn:disc:scalarEMtheorem:Condsb}) 
as a function of $\mu$ for $a=0$ and $b=1$.
} 
\end{figure}

Some numerical results with this estimate are
shown in Fig.~\ref{fig:hatXTypeII}.
We observe that for given $\mu$ and a sufficiently large $s^2$, 
the estimate discretizes, i.e., $\hat x \in \{a, b\}$.
Moreover,
and most importantly 
(and different from estimation as in Section~\ref{sec:disc:JointMAP:Analysis}),
EM converges to $a$ if $\mu$ is closer 
to $a$ than to $b$, and to $b$ if $\mu$ is closer to $b$,
independently of the initialization.
Quantitatively, we have
\begin{theorem}\label{theorem:scalarEM}
Assume $a<b$.
For $\mu < (a+b)/2$, the function
\begin{IEEEeqnarray}{rCl} \label{eqn:disc:scalarEMtheoremIntegral}
\theta \mapsto \int_{-\infty}^\infty p(\breve y, x; \theta) \dd x
= \int_{-\infty}^\infty \Normal{x; \mu, s^2} p(x;\theta) \dd x
\IEEEeqnarraynumspace
\end{IEEEeqnarray}
has a maximum at $\va=0$ and $\vb=(a-b)^2$ (resulting in $\hat x = a$)
and no other extrema if and only if
\begin{IEEEeqnarray}{rCl} \label{eqn:disc:BinCondEMa}
 s^2 > s_{\text{EM}}^2,
\end{IEEEeqnarray}
where 
\begin{IEEEeqnarray}{rCl}\label{eqn:disc:scalarEMtheorem:Condsa}
s_{\text{EM}}^2 \!=\! 
    \left\{ \begin{array}{ll}
       (3\!-\!\sqrt{8})(a\!-\!\mu)(b\!-\!\mu)   & \text{if $\mu \!<\! a \!-\! \frac{|a\!-\!b|}{\sqrt{2}}$}, \\
       \frac{(a-\mu)^2 |a-b|}{(a+b)-2\mu}       & \text{if $a \!-\! \frac{|a\!-\!b|}{\sqrt{2}} \leq \mu \!<\! \frac{a\!+\!b}{2}$}.
   \end{array}\right.  
\end{IEEEeqnarray}
Likewise, for $\mu > (a+b)/2$, (\ref{eqn:disc:scalarEMtheoremIntegral})
has a maximum at $\vb=0$ and $\va=(a-b)^2$ (resulting in $\hat x = b$)
and no other extrema if and only if 
\begin{IEEEeqnarray}{rCl} \label{eqn:disc:BinCondEMb}
 s^2 > s_{\text{EM}}^2, 
\end{IEEEeqnarray}
where 
\begin{IEEEeqnarray}{rCl}\label{eqn:disc:scalarEMtheorem:Condsb}
s_{\text{EM}}^2 \!=\! 
   \left\{ \begin{array}{ll}
       (3\!-\!\sqrt{8})(a\!-\!\mu)(b\!-\!\mu)           & \text{if $\mu \!>\! b \!+\! \frac{|a\!-\!b|}{\sqrt{2}}$}, \\
       \frac{(b-\mu)^2 |a-b|}{2\mu-(a+b)}   & \text{if $\frac{a\!+\!b}{2} \!<\! \mu \leq  b \!+\! \frac{|a\!-\!b|}{\sqrt{2}}$}. 
   \end{array} \right.
    \IEEEeqnarraynumspace
\end{IEEEeqnarray}
\end{theorem}
The proof is not easy%
\footnote{To the reviewers: have a look at the proof, just to get an impression.} 
and does not fit into this paper, 
but can be found in~\cite[App.~C]{keusch2021binaryNUVext} and \cite[App.~B.2]{Keusch2022PhD}. 
The value of $s_{\text{EM}}^2$  
as a function of $\mu$ is plotted in Fig.~\ref{fig:BoundEM}.
For example, $s_{\text{EM}}^2 = 0.225$ for $\mu=0.3$, $a=0$ and $b=1$ (cf. Fig.~\ref{fig:hatXTypeII}).

\subsection{$M$-Level Prior}
\label{sec:MLevel}

\begin{figure}
\centering
\includestandalone[]{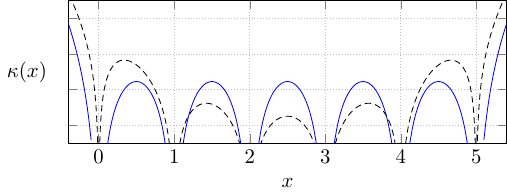}
\caption{\label{fig:effectivePriorTrivialVsComposite}%
Generalization of (\ref{eqn:disc:CostFunctionBinPrior}) to $M=6$ equidistant levels.
Dashed: using (\ref{eqn:disc:trivialPriorModel}).
Solid blue: using (\ref{eqn:disc:SumOfBinaries}) and (\ref{eqn:disc:SumOfBinariesEqualCoeffs}).
}
\end{figure}

An obvious attempt to generalize (\ref{eqn:disc:TwoLvlPrior})
to more than two levels is 
\begin{equation} \label{eqn:disc:trivialPriorModel}
p(x; \theta) 
\eqdef \Normal{x; a, \sigma_a^2} \Normal{x; b, \sigma_b^2} 
    \Normal{x; c, \sigma_c^2} \cdots
\end{equation}
with $\theta \eqdef (\sigma_a^2, \sigma_b^2,\, \ldots)$.
However, this turns out not to work very well
since it introduces a bias towards the levels in the middle range
as illustrated in Fig.~\ref{fig:effectivePriorTrivialVsComposite}.

\begin{figure}
\centering
\includestandalone[]{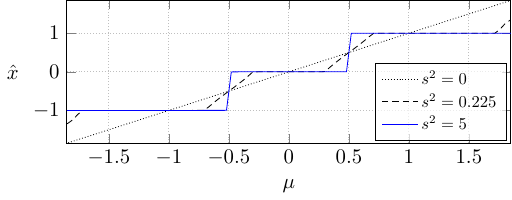}  
\caption{\label{fig:hatXTypeIIThreeLvl}%
Generalization of Fig.~\ref{fig:hatXTypeII} to three levels $\{ -1, 0, 1 \}$
using (\ref{eqn:disc:TernarybySumOfTwo}).}
\end{figure}

Good results are obtained with linear combinations of auxiliary binary (or binarized)
variables. For example, 
constraining $X$ to three levels $\{ -b, 0, b \}$
can be achieved by writing 
\begin{equation} \label{eqn:disc:TernarybySumOfTwo}
X = b X_1 - b X_2
\end{equation}
where both $X_1$ and $X_2$ are constrained to $\{ 0, 1 \}$ 
by means of independent priors (\ref{eqn:disc:TwoLvlPrior}), i.e.,
\begin{IEEEeqnarray}{rCl}
p(x_1, x_2; \theta_1, \theta_2) 
 & = & 
  \Normal{x_1; 0, \sigma_{1,a}^2} \Normal{x_1; 1, \sigma_{1,b}^2} \nonumber \\
  && {}\cdot\Normal{x_2; 0, \sigma_{2,a}^2} \Normal{x_2; 1, \sigma_{2,b}^2}.
     \IEEEeqnarraynumspace
\end{IEEEeqnarray}
The corresponding generalization of Fig.~\ref{fig:hatXTypeII}
is shown as solid line in Fig.~\ref{fig:hatXTypeIIThreeLvl}.

More generally, we can write $X$ as a linear combination
\begin{equation} \label{eqn:disc:SumOfBinaries}
X = \sum_{j=1}^J \beta_j X_j + \beta_0
\end{equation}
of independent binary (i.e., binarized to $\{ 0, 1\}$) variables $X_1,\, \ldots, X_J$.
The choice of $J$ and of the coefficients $\beta_0,\, \ldots, \beta_J$
is highly nonunique. 
Choosing $\beta_j = 2^{j-1}$ for $j>0$ does not work well empirically.
Good results are obtained with 
\begin{equation} \label{eqn:disc:SumOfBinariesEqualCoeffs}
\beta_1 = \ldots = \beta_J,
\end{equation}
resulting in $M=J+1$ equidistant levels for $X$.
(Related representations were used in \cite{FrLg:sradda2006}.)
The corresponding generalization of (\ref{eqn:disc:CostFunctionBinPrior})
is illustrated in Fig.~\ref{fig:effectivePriorTrivialVsComposite}.

In (\ref{eqn:disc:SumOfBinaries}), 
\mbox{$X_1=0$} and \mbox{$X_2=1$} has the same effect on $X$ as $X_1=1$ and $X_2=0$. 
The estimation algorithm must somehow 
choose
among such equivalent configurations. 
However, depending on the details of the implementation,
the estimation algorithm may not, by itself, be able to break such symmetries.
This problem can be solved
by a slightly asymmetric initialization of the variances,
e.g., 
\begin{equation}
\sigma_{1,a}^2 = \sigma_{1,b}^2 \neq \sigma_{2,a}^2 = \sigma_{2,b}^2,
\end{equation}
where the inequality is almost an equality.

\section{Application Examples}
\label{sec:app}

We now demonstrate the versatility of the proposed approach 
(as described in Section~\ref{sec:ProposedApproach})
by sketching its application to some exemplary control problems.
For an in-depth study of an industrial control problem,
the reader is referred to the companion paper~\cite{keusch2023ContrSysTechn}.


\subsection{Squared-Error Fitting with Binary Control} \label{sec:DAC}

We begin with an example like Example~\ref{example:prob:binaryControl}
of Section~\ref{sec:SysApproachExamples}:
we wish to steer a linear system as in Section~\ref{sec:SysApproachExamples} 
with a $\{ 0, 1\}$-valued control signal 
$u_1,\ldots,u_K$ such that its scalar output $y_1,\ldots,y_K$
follows a given target $\breve y_1, \ldots, \breve y_K$ 
such that 
\begin{equation} \label{ex:DAC:fittingCost}
\sum_{k=1}^K (y_k - \breve y_k)^2
\end{equation}
is as small as possible.
We do not actually aim for the global minimum of (\ref{ex:DAC:fittingCost}), 
but we hope to get close to it.

The quadratic penalty~(\ref{ex:DAC:fittingCost}) is readily 
expressed by $\msgb{p}{}(y_k; \theta_{Y_k})$ as in (\ref{eqn:p:msgbYk})
with fixed parameters 
$\msgb{m}{Y_k} = \breve y_k$ and $\msgb{V}{Y_k} = \msgb{\sigma}{Y}^2 > 0$. 
The choice of $\msgb{\sigma}{Y}^2$ will be discussed below.

The constraint $u_k\in \{0,1\}$ is
expressed by $\msgf{p}{}(u_k;\theta_{U_k})$ as in (\ref{eqn:p:msgfUk}) 
with unknown parameters~$\msgf{m}{U_k}$ and~$\msgf{V}{U_k}$.
In Step~2 of IAKE, $\msgf{m}{U_k}$ and~$\msgf{V}{U_k}$
are updated 
using (\ref{eqn:TableNUP:Bin})
with $X=U_k$, $a=0$, $b=1$, 
$m_X = m_{U_{k}}^{(i)}$, and $V_X = V_{U_k}^{(i)}$.

For the numerical experiments, we use a stable linear system 
(or linear filter)
with transfer function
(=~the Laplace transform of the impulse response)
\begin{IEEEeqnarray}{rCl} \label{eqn:dac:transferFunction}
  G(s) = \frac{35037.9}{s^3 + 71.9 s^2 + 2324.8s + 35037.9}
\end{IEEEeqnarray}


The transfer function~(\ref{eqn:dac:transferFunction}) is transformed into 
state-space form and discretized 
using a sampling interval of $T = 0.003$ seconds, resulting in a 
discrete-time system as in~(\ref{eqn:lssm:detLSSM}) with 
state space dimension $N=3$ and matrices
\begin{IEEEeqnarray}{rCl} \IEEEyesnumber \IEEEyessubnumber*
    A \!&=&\! \bma
            0.7967&  -6.3978& -94.2123\\
            0.0027 &  0.9902 & -0.1467\\
            0      &  0.0030 &  0.9999
          \ema, \; 
  B \!=\! \bma 0.0027 \\ 0 \\ 0  \ema,  \IEEEeqnarraynumspace
\end{IEEEeqnarray}
  and
\begin{IEEEeqnarray}{rCl}
  \IEEEyessubnumber*
  C &=& \bma 0 & 0 & 35037.9 \ema.
\end{IEEEeqnarray}

\begin{figure}
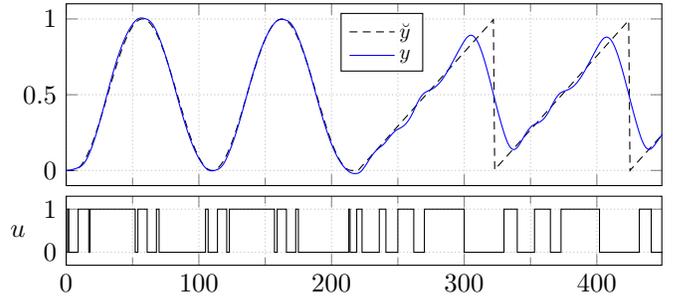

  \centering
  \includestandalone{inc_pdf/dac}
  \caption{
  Binary-input control (or digital-to-analog conversion) as in Section~\ref{sec:DAC}
  with target waveform $\breve y$ (dashed),
  binary control signal $u$ computed by the proposed algorithm (bottom), 
  and resulting output signal y (solid blue).
  }
  \label{fig:prob:dac}
\end{figure}

The numerical results shown in Fig.~\ref{fig:prob:dac} are obtained with 
$\msgb{V}{Y_k} = \msgb{\sigma}{Y}^2 = 0.045$ and $K=450$.  
In the first half of Fig.~\ref{fig:prob:dac} the target waveform 
can be well approximated; 
in the second half of Fig.~\ref{fig:prob:dac}, 
the target waveform falls outside the passband 
of the filter (\ref{eqn:dac:transferFunction}).

Constraint satisfaction 
can be controlled by $\msgb{\sigma}{Y}^2$: 
if the final estimate $u_k$ fails to satisfy $u_k\in\{0,1\}$ for all~$k$, 
$\msgb{\sigma}{Y}^2$ should be increased, e.g., by a factor of~2,
cf.\ Section~\ref{sec:ConstraintsChecking}.
(But $\msgb{\sigma}{Y}^2$ should not be 
chosen to be unnecessarily large since this slows down the
convergence of IAKE.)

\begin{figure}
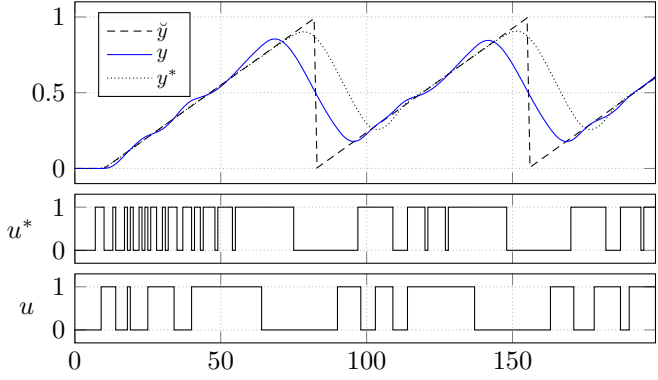

\begin{center}
\includestandalone{inc_pdf/3th_ord_SK_filter_vs_optimal_with_w_8}
\caption{\label{fig:versus_optimal}%
Comparing the proposed method (with planning horizon $K=200$) 
with an optimal (exhaustive search) controller
with planning horizon $K=8$.  
The former yields a significantly better approximation 
(solid blue $y$ with $\text{MSE}=0.01972$) than the latter  
(dotted $y^*$ with $\text{MSE}=0.04885$).
}
\end{center} 
\end{figure}

\begin{figure}
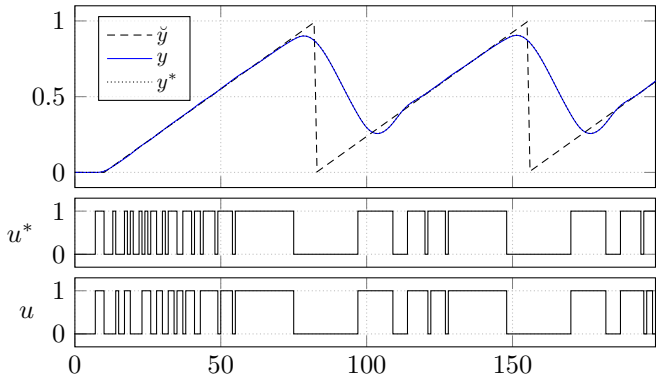

\begin{center}
\includestandalone{inc_pdf/3th_ord_SK_filter_with_w_8_vs_optimal_with_w_8}
\caption{\label{fig:versus_optimal_both_short_horizon}%
Comparing the proposed method (solid blue $y$)
with an optimal (exhaustive search) controller (dotted $y^\ast$, covered by $y$),
both in receding-horizon mode with planning horizon $K=8$.
The approximation error is nearly identical  
($\text{MSE}=0.04899$ vs.\ $\text{MSE}=0.04885$).
}
\end{center} 
\end{figure}

A comparison with an ``optimal'' controller is shown in 
Figs.\ \ref{fig:versus_optimal} and~\ref{fig:versus_optimal_both_short_horizon}.
This ``optimal'' controller determines the binary control sequence 
by exhaustive search, which severely limits its planning horizon $K$.
By contrast, the proposed method can work with a full-length planning horizon.
Fig.~\ref{fig:versus_optimal} illustrates the advantage of the latter.
But even if both methods work with the same (short) planning horizon,
Fig.~\ref{fig:versus_optimal_both_short_horizon} shows that,
in this example, 
the proposed method yields an essentially optimal control input.

\subsection{Corridor Control with Different Input Constraints}
\label{sec:var:boxOutput}

Assume we wish to keep the system output $y$ within a corridor around a target 
$\breve y$, i.e., 
we wish $y_k$ to satisfy
\begin{IEEEeqnarray}{rCl}
  a_k \leq y_k - \breve y_k \leq b_k, \quad k \in \{1, \dots, K\},
\end{IEEEeqnarray}
for fixed bounds \mbox{$a_k, b_k \in \Reals$}.
These constraints can be expressed 
by $\msgb{p}{}(y_k;\theta_{Y_k})$ as in (\ref{eqn:p:msgbYk})
with unknown parameters $\msgb{m}{Y_k}$ and $\msgb{V}{Y_k}$
that are updated (in Step~2 of IAKE) 
using (\ref{eqn:TableNUP:box})
with $X=Y_k$, 
$\msgf{m}{X}=\msgb{m}{Y_k}$, $\msgf{V}{X}=\msgb{V}{Y_k}$, 
$a = a_k + \breve y_k$, $b = b_k + \breve y_k$, 
$m_X = m_{Y_{k}}^{(i)}$, $V_X = V_{Y_k}^{(i)}$,
and slope parameter $\gamma = \gamma_{Y}$.

In the following, we consider five different version of this problem,
with different constraints on the control signal $u$,
as illustrated in Fig.~\ref{fig:var:boxOutput}
(with numerical values given below).
Note that (\ref{enumi:Corridor:sparseInput}), (\ref{enumi:Corridor:ternaryInput}),
and (\ref{enumi:Corridor:sparseInputChanges})
amount to nonconvex optimization problems.

{\renewcommand{\theenumi}{V\arabic{enumi}}
\begin{enumerate}
\item\label{enumi:Corridor:L2Input}
The input $u$ is regularized by an $L_2$ penalty, which is expressed
by $\msgf{p}{}(u_k;\theta_{U_k})$ with fixed parameters $\msgf{m}{U_k}=0$
and $\msgf{V}{U_k} = \sigma_U^2$.

Constraint satisfaction can be guaranteed by increasing either $\gamma_Y$
or $\sigma_U^2$ (e.g., by a factor of~2), if necessary, 
as described in Section~\ref{sec:ConstraintsChecking}.
Since the optimization problem is convex, 
the choice of these parameters has no effect on the estimate $u$ 
as long as the constraints are satisfied. 
However, choosing these parameters unnecessarily large
makes the convergence of IAKE unnecessarily slow.

\item\label{enumi:Corridor:InputSlopeGeq}
The input $u$ is required to satisfy 
\begin{equation}\label{eqn:Corridor:InputDerivGeq}
\tilde u_k \eqdef u_k - u_{k-1} \geq a
\end{equation}
for all~$k$.
%
To this end,
we modify the state space model (\ref{eqn:lssm:detLSSM})
to 
$\tilde x \eqdef \bma u_k, x_k \ema^\T$, 
new input $\tilde u_k$,
and
\begin{IEEEeqnarray}{rCl} \label{eqn:Corridor:AugmentedLSSM}
  \tilde A = \bma 1 & 0_{1 \times N} \\ B & A \ema, \; 
  \tilde B = \bma 1 \\ 0_{N \times 1} \ema, \; 
  \tilde C = \bma 0 & C \ema. \IEEEeqnarraynumspace
\end{IEEEeqnarray}
(But the actual control signal is still $u_k$, i.e., the first component of $\tilde x_k$.)

The constraint (\ref{eqn:Corridor:InputDerivGeq}) is then expressed by 
$\msgf{p}{}(\tilde u_k;\theta_{\tilde U_k})$ as in (\ref{eqn:p:msgfUk}),
where $\msgf{m}{\tilde U_k}$ and~$\msgf{V}{\tilde U_k}$ are updated using (\ref{eqn:TableNUP:geq}).

If the constraints can be satisfied at all,
then IAKE (with sufficiently large $\gamma_Y$ and $\gamma_{\tilde U}$, cf.\ Section~\ref{sec:ConstraintsChecking})
will find a pertinent control signal $u$.

\item\label{enumi:Corridor:sparseInput}
The input $u$ is required to be sparse,
which is achieved by $\msgf{p}{}(u_k;\theta_{U_k})$ as in (\ref{eqn:p:msgfUk})
with parameters~$\msgf{m}{U_k}=0$ and $\msgf{V}{U_k}$ updated by~(\ref{eqn:TableNUV:plainNUV})
(with $\msgf{V}{U_k} = \msgf{V}{X}$, $V_X = V_{U_k}$, and $m_X = m_{U_k}$).
Constraint satisfaction is enforced by sufficiently large $\gamma_Y$.

A variation of this problem was discussed in~\cite{Hoffmann2017}.
Empirically, (\ref{eqn:TableNUV:plainNUV}) works better than 
standard $L_1$ regularization~\cite{Tibshirani1996}, which is effected by (\ref{eqn:TableNUV:L1}).

\item\label{enumi:Corridor:ternaryInput}
The input $u$ is required to satisfy $u_k \in \{-1, 0, 1\}$ 
for all~$k$.
This is achieved with $u_k = \tilde u_{k}^+ -  \tilde u_{k}^-$,
where both $\tilde u_{k}^+$ and $\tilde u_{k}^-$ are constrained 
to $\{ 0, 1 \}$ by (\ref{eqn:TableNUP:Bin}),
as described in Section~\ref{sec:MLevel}.

Constraint satisfaction can be encouraged with sufficiently large $\gamma_Y$,
but cannot actually be guaranteed (even if the problem itself is feasible).
However, empirically, the proposed method works very well also in this case. 

\item\label{enumi:Corridor:sparseInputChanges}
$\tilde u_k \eqdef u_k - u_{k-1}$ is required to be sparse. 
This is achieved with the modified state space model (\ref{eqn:Corridor:AugmentedLSSM})
and a sparsifying prior on $\tilde U_k$ as in (\ref{enumi:Corridor:sparseInput}).
\end{enumerate}
}

The numerical results in Fig.~\ref{fig:var:boxOutput} 
are obtained with the state space model 
\begin{IEEEeqnarray}{rCl} \IEEEyesnumber \IEEEyessubnumber*
  A &=& \bma 1  & 0 & 0 \\ 1 & 1 & 0 \\ 1/2 & 1 & 1 \ema, \quad
  B = 0.0015 \bma 1 \\ 1/2 \\ 1/3 \ema,
\end{IEEEeqnarray}
and 
\begin{IEEEeqnarray}{rCl} \IEEEyessubnumber*
  C = \bma 0 & 0 & 1 \ema.
\end{IEEEeqnarray}
Furthermore, we have $K=175$, 
$\gamma_Y = 10$ throughout,
$\sigma_U^2 = 10$
in (\ref{enumi:Corridor:L2Input}),  
and $a = -0.03$ 
and $\gamma_U = 10$ 
in (\ref{enumi:Corridor:InputSlopeGeq}). 

\begin{figure}
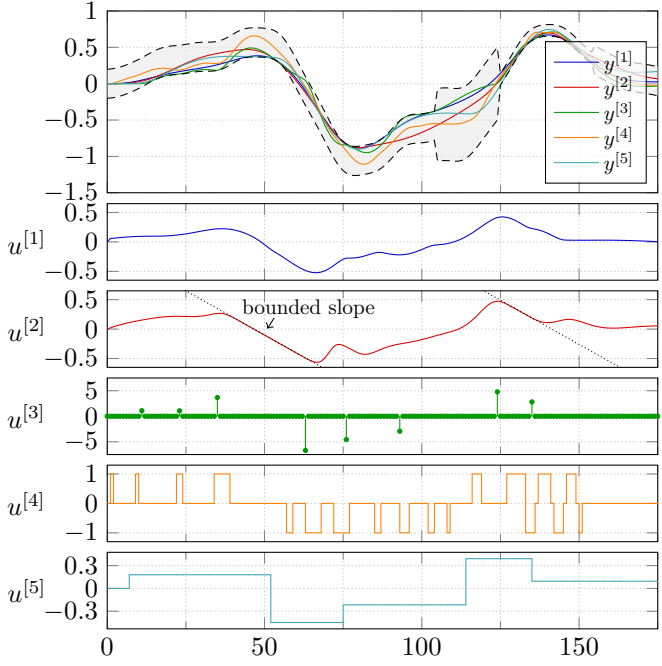

  \centering
  \includestandalone{inc_pdf/boxConstraintOutput}
  \caption{\label{fig:var:boxOutput}%
  Corridor control with different constraints on the input.
  Top row: prescribed corridor (dashed) and resulting output signals.
  Other rows, from top to bottom:
  input $u$ with $L_2$ penalty (\ref{enumi:Corridor:L2Input}),
  input $u$ with lower-bounded slope (\ref{enumi:Corridor:InputSlopeGeq}),
  input $u$ with sparsifying penalty (\ref{enumi:Corridor:sparseInput}),
  ternary input $u$ (\ref{enumi:Corridor:ternaryInput}),
  input $u$ with sparsifying penalty on level switches (\ref{enumi:Corridor:sparseInputChanges}).
  }
\end{figure}
\noindent


\subsection{Double-Slit Flappy Bird Control}
\label{sec:flappyBird}

The following control problem is a variation of the 
\emph{flappy bird} computer game~\cite{wiki_flappy_bird}.
(This example improves on the related example in~\cite{keusch2021binaryNUV}, 
which did not use box constraints.)

\begin{figure}
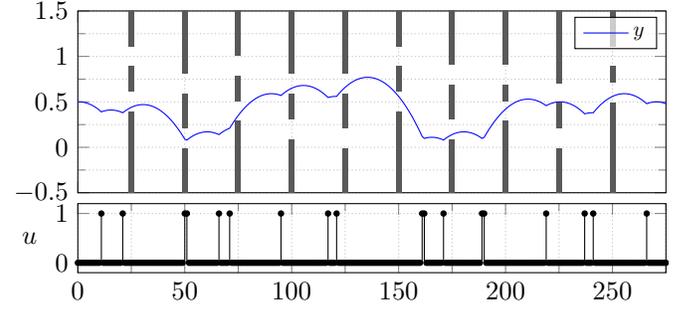

  \centering
  \includestandalone{inc_pdf/flappy_bird_two_slits}
  \caption{\label{fig:traj:FlappyBirdDoubleSlits}%
  Double-slit flappy bird control with binary control signal $u$ and
  resulting trajectory $y$. 
  }  
\end{figure}

Consider a physical system consisting 
of a point mass $m$ moving forward (left to right in Fig.~\ref{fig:traj:FlappyBirdDoubleSlits}) 
with constant horizontal velocity
and ``falling'' vertically with constant acceleration~$g$. 
The $\{0,1\}$-valued control signal $u$ affects the system only if $u_k=1$,
in which case a fixed value is added to the vertical momentum.
We wish to steer the point mass such that it passes 
through a sequence of double slits
as illustrated in Fig.~\ref{fig:traj:FlappyBirdDoubleSlits}.

For this example, we need a slight generalization
of~(\ref{eqn:lssm:detLSSM}) as follows.
The state $x_k \in\Reals^2$ (comprising the vertical position and the vertical velocity) 
evolves according to
\begin{IEEEeqnarray}{rCl} 
x_k & = & \bma 1 & T \\ 0 & 1 \ema x_{k-1}
          + \bma 0 \\ 1/m \ema u_k + \bma 0 \\ -Tg \ema.
\end{IEEEeqnarray}
The output $y_k \eqdef \bma 1 & 0 \ema x_k$ is the vertical position.

However, we directly constrain not $y_k$, but an auxiliary output $\tilde y_k$
that is defined as follows. 
Let $\mathcal S \subset \{ 1, \ldots, K \}$ be the positions of the double slits.
For $k \not\in \mathcal S$, there is no output $\tilde y_k$;
for $k \in \mathcal S$,
\begin{IEEEeqnarray}{rCl} \label{eqn:FlappyShiftAug}
  \tilde y_k \eqdef y_k + s_k,
\end{IEEEeqnarray}
where $s_k \in \{0, d_k\}$ ``selects'' either the lower or
the upper slit, and where ${d_k \in \Reals}$ specifies the vertical distance between them.

The double-slit constraint
\begin{IEEEeqnarray}{rCl} \label{eqn:traj:doubleSlitConstraint}
  y_k \in [a_k, b_k]  \quad \text{or} \quad y_k \in [a_k-d_k, b_k-d_k]
\end{IEEEeqnarray}
is then expressed by a box constraint on $\tilde Y_k$ (with bounds $a_k$ and $b_k$)
and a $\{ 0, d_k \}$-constraint on $S_k$.

The constraints on the input are simply expressed with 
a binarizing prior on $U_k$ with levels $\{0, 1\}$ for all $k$.

The numerical results in Fig.~\ref{fig:traj:FlappyBirdDoubleSlits} are obtained with 
$K = 300$, $m = 1$, $T = 0.1$, $g = 0.2$, $\gamma = 100$, and $a_k, b_k$ and 
$d_k$ according to Fig.~\ref{fig:traj:FlappyBirdDoubleSlits}.

\subsection{Trajectory Planning with Obstacle Avoidance}

Consider the following situation. 
An object is moving in a two-dimensional plane.
Its position at time $k$ is $y_k= \bma y_{k,1}, y_{k,2} \ema^\T \in \Reals^2$,
which is governed by the
state space model~(\ref{eqn:lssm:detLSSM}) with
\begin{IEEEeqnarray}{rCl} 
  A &=& \bma  1 & 0 & 0 & 0 \\
              T & 1 & 0 & 0 \\
              0 & 0 & 1 & 0 \\
              0 & 0 & T & 1 
       \ema, \;
  B = \bma    T & 0 \\
              0 & 0 \\
              0 & T \\
              0 & 0
      \ema, \;
  C = \bma 0 \\ 1 \\ 0 \\ 1 \ema^\T \!\!, \IEEEeqnarraynumspace
\end{IEEEeqnarray}
where $T$ is the discretization interval
and $u_k = \bma u_{k,1}, u_{k,2} \ema^\T$ is the acceleration.

\begin{figure}
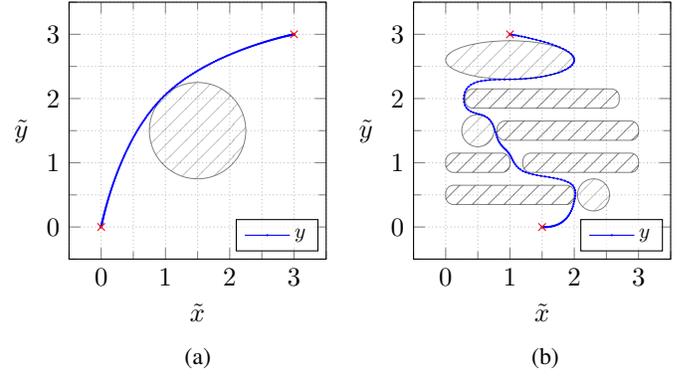
   
  \hspace*{-1.7em}
  \begin{subfigure}[b]{0.49\linewidth}
    \includestandalone{inc_pdf/keepOutZone_spheres_single}
    \caption{\hspace*{-5em}}\label{fig:traj:koz_spheres_single} 
  \end{subfigure}
\begin{subfigure}[b]{0.49\linewidth}
    \includestandalone{inc_pdf/keepOutZone_mixed_I}
    \caption{\hspace*{-6em}}\label{fig:traj:koz_mixed_I} 
\end{subfigure}
\caption{(a) Trajectory planning with a single spherical obstacle at 
$\bma 1.5, 1.5 \ema^\T$. 
Note that the optimal trajectory $y$ is not unique. (b) Trajectory planning with obstacles of various shapes. 
The obtained trajectory $y$ is only locally optimal.}
\end{figure}

Assume we wish to plan a trajectory starting from $\bma 0, 0 \ema^\T$ (with zero velocity) 
and ending at  $ \bma 3,3 \ema^\T$ (with zero velocity), 
while avoiding a spherical obstacle 
at $c = \bma 1.5, 1.5 \ema^\T$ with radius $r = 0.75$ (see Fig.~\ref{fig:traj:koz_spheres_single}).
In addition, we wish to minimize the squared norm of the acceleration, i.e.,
\begin{IEEEeqnarray}{rCl}
  \sum_{k=1}^K \| u_k\|^2,
\end{IEEEeqnarray}
which is easily handled by a zero-mean Gaussian prior on $U_k$, for all $k$. 

\begin{figure}
  \centering
  \includestandalone[scale=0.75]{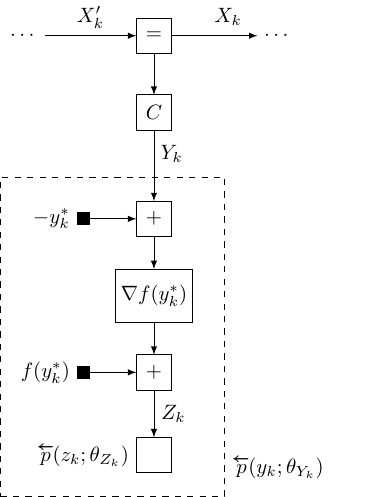}
  \caption{\label{fig:traj:FGkozSingleObstacle}%
  Factor graph of the half-space prior $\msgb{p}{}(z_k; \theta_{Z_k})$
  on the linearized observation~(\ref{eqn:ExObstacles:LinearizedObs}).
  }
\end{figure}

The obstacle can be avoided by a half-space constraint on the auxiliary variable 
\begin{IEEEeqnarray}{rCl} \label{eqn:traj:nonLinFunc}
  \tilde z_k \eqdef  \| y_k - c\| = f(y_k), \quad k \in \{1, \dots, K\},
\end{IEEEeqnarray}
which is the distance from $y_k$ to the center $c$ of the obstacle. 
Specifically, we use a half-space NUV prior to enforce
\begin{IEEEeqnarray}{rCl}
  \tilde z_k > r.
\end{IEEEeqnarray}
It remains to deal with the problem that~(\ref{eqn:traj:nonLinFunc}) 
is a nonlinear function of $y_k$.
We solve this problem in the most obvious way,
by using the linearization 
\begin{IEEEeqnarray}{rCl}  \label{eqn:ExObstacles:LinearizedObs}
  z_k & = & f(y_k^*) + \nabla f(y_k^*) (y_k - y_k^*) \approx f(y_k)
\end{IEEEeqnarray}
(as illustrated in Fig.~\ref{fig:traj:FGkozSingleObstacle}),
where $y_k^* \in \Reals^2$ is the previous estimate of $Y_k$ and $\nabla f(y_k^*)$ 
is the gradient of $f$ at $y_k=y_k^*$.
%

The numerical results illustrated in Fig.~\ref{fig:traj:koz_spheres_single} 
are obtained 
with $T=1$, $\gamma = 5$, $\msgf{m}{U_k} = \bma 0, 0 \ema^\T$, 
$\msgf{V}{U_k} = \diag{0.1, 0.1}$, 
and boundary conditions 
\begin{IEEEeqnarray}{rCl}
  \msgf{m}{X_0} &=& \bma 0 , 0 , 0 , 0 \ema^\T, \\
  \msgb{m}{X_K} &=& \bma 0 , 3 , 0 , 3 \ema^\T, \quad \text{and} \\
  \msgf{V}{X_0} &=& \msgb{V}{X_K} = 0_{4 \times 4}.
\end{IEEEeqnarray}
Note that the optimal solution of the given problem is not unique since the problem is 
geometrically symmetric. The obtained solution depends on the initial conditions. 

The method of this example is easily extended to multiple obstacles 
by concatenating multiple instances of the part shown 
in Fig.~\ref{fig:traj:FGkozSingleObstacle}.
The method is not limited to spherical obstacles as long as the nonlinearity of $f$ is good-natured. 
Ellipses, squares, rectangles, and linear transformations (e.g., scaling and rotations) 
thereof have been successfully implemented 
by choosing $f$ accordingly. An example with multiple obstacles of various shapes is given in 
Fig.~\ref{fig:traj:koz_mixed_I}, the details are omitted.


\subsection{Minimal-Time Race Track Control} \label{sec:racetrack}

Autonomous racing is a version of autonomous driving where the goal is to 
complete a given race track in the shortest time possible.
The following challenges must be dealt with:
\begin{itemize}
  \item Nonlinear vehicle dynamics.
  \item Physical limitations of the vehicle such as maximal steering angle and maximal motor torque.
  \item Collision avoidance with track boundaries.
\end{itemize}
Several methods to solve this control problem have been proposed in the 
literature~\cite{Qian2016, Qian2016a, Rosolia2019}. 
We now show how this problem can be addressed with the approach of this paper.

\subsubsection{State Space Models: Step by Step}

As recommended in the literature,
we will use a curvilinear coordinate 
system~\cite{Lot2014, Micaelli1993, Lenain2007, Lenain2008}, which simplifies 
expressing the constraints imposed by the track boundaries.

We begin by describing the vehicle dynamics using the standard 
\textit{Ackermann vehicle model}~\cite{rajamani2011vehicle}
in Cartesian coordinates 
$(\tilde x, \tilde y)\in\Reals^2$,
from which the final state space model will be obtained in a series of transformations.
We thus begin with
the differential equation
\begin{IEEEeqnarray}{rCl} \label{eqn:race:ackermannDynamicsCartesian}
  \frac{\dd x}{\dd t} & = &  f(x(t), u(t)) = 
  \bma 
    v(t) \cos(\theta(t)) \\ 
    v(t) \sin(\theta(t)) \\ 
    v(t) \frac{\tan(\delta(t))}{\ell} \\ 
    a(t) \\
    \dot a(t) \\
    \dot \delta(t) 
  \ema
  \IEEEeqnarraynumspace
\end{IEEEeqnarray}
with state
\begin{IEEEeqnarray}{rCl}
  x(t) & = & \bma \tilde x(t) , \tilde y(t) , \theta(t) , v(t) , a(t) , \delta(t) \ema^\T,
  \IEEEeqnarraynumspace
\end{IEEEeqnarray}
input
\begin{IEEEeqnarray}{rCl}
  u(t) & = & \bma \dot \delta(t) \\ \dot a(t) \ema,
      \label{eqn:RaceTrack:DerivativesInput}
\end{IEEEeqnarray}
and 
heading angle~$\theta$, 
(front wheel) steering angle $\delta$, 
vehicle length~$\ell$, 
speed~$v$, and acceleration~$a$,
and where the dot in $\dot a(t)$ etc.\ denotes the derivative with respect to the time~$t$.

Using (\ref{eqn:RaceTrack:DerivativesInput}), rather than $\delta(t)$ and $a(t)$ directly,
as inputs is minor embellishment: it will allows us to discourage
very rapid changes of $\delta(t)$ and $a(t)$ by a suitable penalty.

\begin{figure}
  \center
  \includestandalone[scale=0.8]{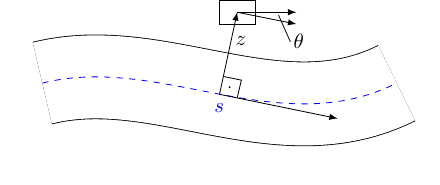}
  \caption{Curvilinear coordinate system, where $s$ is the progress along the 
  center line, $z$ is the vehicle's displacement perpendicular to the centerline, 
  and $\theta$ is the heading angle relative to a tangent vector at $s$.}
  \label{fig:race:curvilinear}
\end{figure}

In a next step, we transform this state space 
model into a curvilinear coordinate system
as illustrated in Fig.~\ref{fig:race:curvilinear}.
%
The first coordinate $s$ (of this curvilinear coordinate system) 
is the progress along the center line of the race track.
The second coordinate $z$ is the perpendicular distance of the
vehicle to the center line at~$s$. 
The (new) angle $\theta$ is the angle between 
the vehicle's direction of travel and the tangent vector at~$s$.
Consequently, the (new) state vector is
\begin{IEEEeqnarray}{rCl}
  x(t) = 
  \bma 
    s(t) , 
    z(t) , 
    \theta(t) , 
    v(t) ,
    a(t) , 
    \delta(t)  
  \ema^\T.
\end{IEEEeqnarray}
In these new coordinates, 
the vehicle dynamics~(\ref{eqn:race:ackermannDynamicsCartesian}) 
are given by
\begin{IEEEeqnarray}{rCl} \label{eqn:race:ackermannDynamicsCurvilinearTime}
  \frac{\dd x}{\dd t}
  = f_t(x(t), u(t)) 
  = \bma  
      v \frac{\cos(\theta)}{\imk} \\
      v \sin(\theta) \\
      v \left ( \frac{\tan(\delta)}{\ell} - \frac{\kappa(s) \cos(\theta)}{\imk} \right) \\
      a \\
      \dot a \\ 
      \dot \delta
    \ema, \IEEEeqnarraynumspace
\end{IEEEeqnarray}
where $\kappa(s)$ is the curvature of the center line.
Note that the right-hand side of~(\ref{eqn:race:ackermannDynamicsCurvilinearTime}) 
depends on the time $t$, which is omitted for readability.

In~(\ref{eqn:race:ackermannDynamicsCurvilinearTime}), the
independent variable is time, 
which is inconvenient for minimal-time optimization.
We therefore transform the state space model once more, 
into a form where the independent variable is $s$.
The transformed model follows directly from
\begin{IEEEeqnarray}{rCl}
  \frac{\dd x}{\dd s} 
  &=& \frac{\dd x}{\dd t} \frac{\dd t}{\dd s} 
  =  \left (\frac{\dd s}{\dd t} \right)^{-1} \frac{\dd x}{\dd t} \\
  &=& \left ( \frac{\imk}{v \cos(\theta)} \right) f_t(x(s), u(s)).
\end{IEEEeqnarray}
Accordingly, the new state and input vectors are no longer functions of $t$, but functions 
of $s$, i.e., $x(s)$, and $u(s)$, respectively.
Since $s$ is now the independent variable, we drop the first state and add time
as an additional state, i.e.,
\begin{IEEEeqnarray}{rCl}
  x(s) = 
  \bma 
    z(s) , 
    \theta(s) , 
    v(s) ,
    a(s) ,
    \delta(s) , 
    t(s) 
  \ema^\T.
\end{IEEEeqnarray}
The new model dynamics are
\begin{IEEEeqnarray}{rCl}
  \frac{\dd x}{\dd s}\!=\! f_s(x(s), u(s)) \!=\! 
  \frac{1 \!-\! \kappa(s) z }{v \cos(\theta)} \!
  \bma 
  v \sin(\theta) \\ 
  v \! \left (\! \frac{ \tan(\delta)}{\ell} \!-\! \frac{\cos(\theta)}{\kappa(s)^{-1} \!-\! z} \! \right ) \!\! \\
  a \\
  \dot a  \\
  \dot \delta \\
  1
\ema \!. \IEEEeqnarraynumspace
  \label{eqn:race:nonLinModel}
\end{IEEEeqnarray}
In order to impose suitable state constraints, we define 
a system output
\begin{IEEEeqnarray}{rCl} \label{eqn:race:outputModel}
  y(s) = f_o(x(s)) = 
  \bma 
  z \\ 
  a \\ 
  \delta \\ 
  a^2 + \psi \frac{v^4}{\ell^2} \tan(\delta)^2 
  \ema, 
\end{IEEEeqnarray}
where the last component of~(\ref{eqn:race:outputModel}) is the squared total acceleration~$a_{\text{tot}}^2$, 
and where $\psi$ is a weighting factor to incorporate all unmodeled physical properties of the vehicle.

\begin{figure}
  \center
  \includestandalone[scale=1]{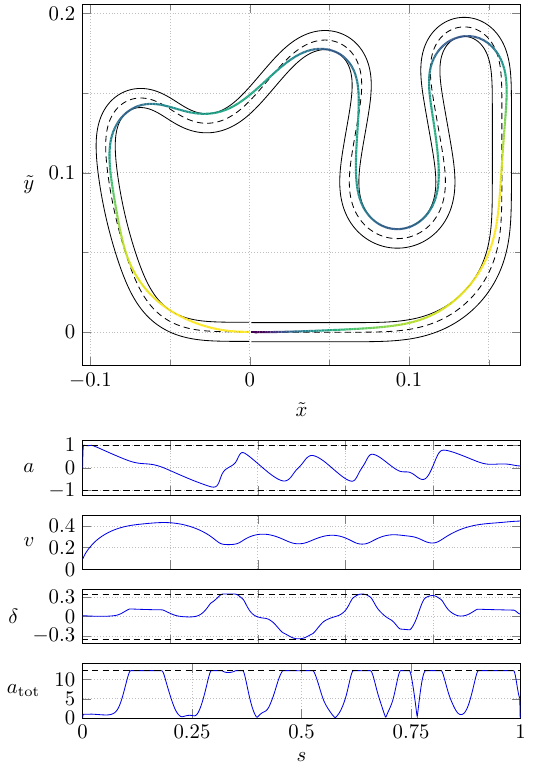}
  \caption{\label{fig:race:raceTrackTotalAcc}%
  Minimal-time racing with constrained longitudinal acceleration $a$, 
  steering angle $\delta$, 
  and total acceleration $a_{\text{tot}}$. 
  The color of the resulting trajectory (top) indicates the speed~$v$.}
\end{figure}

In a final step, we linearize the nonlinear 
model~(\ref{eqn:race:nonLinModel}) and~(\ref{eqn:race:outputModel}) around 
the linearization point $(x^*, u^*)$, yielding the model 
\begin{IEEEeqnarray}{rCl}  \IEEEyesnumber \phantomsection \label{eqn:race:linearModel}  \IEEEyessubnumber*
  \frac{\dd x}{\dd s} &=& \tilde A (x(s) - x^*) + \tilde B(u(s) - u^*) + f_s(x^*, u^*) \IEEEeqnarraynumspace \\
  y(s) &=& \tilde C (x(s) - x^*) + f_o(x^*), 
\end{IEEEeqnarray}
with
\begin{IEEEeqnarray}{rCl}
  \tilde A &=& \frac{\partial f_s(x^*, u^*)}{\partial x}, \; 
  \tilde B = \frac{\partial f_s(x^*, u^*)}{\partial u}, \;
  \tilde C = \frac{\partial f_o(x^*)}{\partial x}. \IEEEeqnarraynumspace
\end{IEEEeqnarray}
The linear model~(\ref{eqn:race:linearModel}) is then discretized 
using a first-order approximation (Euler method), resulting in
\begin{IEEEeqnarray}{rCl} \label{eqn:race:lssm}
  x_{k+1} &=& A (x_k \!-\! x_k^*) + B (u_k \!-\! u_k^*) + x_k^* + T_s  f_s(x_k^*, u_k^*) \IEEEeqnarraynumspace \\
  y_{k} &=& C (x_k \!-\! x_k^*) + f_o(x_k^*),
\end{IEEEeqnarray}
with
\begin{IEEEeqnarray}{rCl}
  A = 1 + T_s \tilde A, \quad 
  B = T_s \tilde B, \quad \text{and} \quad 
  C = \tilde C,
\end{IEEEeqnarray}
and where $T_s$ is the spatial sampling interval.

\subsubsection{Adding the Constraints}

Keeping the vehicle within the track boundaries is achieved by imposing 
box constraints on $z_k$ (the discretized version of $z(s)$) along the track.
Further box constraints on the longitudinal acceleration $a_k$ 
and the steering angle $\delta_k$ enforce physical limitations of the vehicle. 
A box constraint on the total acceleration $a_{\text{tot}}^2$
prevents the vehicle from slipping. 
Finally, minimizing the track 
time is handled by imposing a zero-mean Gaussian penalty 
on the time of arrival $x_{K, 6}$ ($=$ the last component of the state $x_k$ at time $k=K$).

\subsubsection{Numerical Example}

The example shown in Fig.~\ref{fig:race:raceTrackTotalAcc} was obtained with 
the following numerical values:
We use box priors on the corresponding model outputs
to constrain the deviation from the centerline 
to $-0.006 \leq z \leq 0.006$
with $\gamma_z = 0.005$, 
the vehicle's longitudinal acceleration
to $-1 \leq a \leq 1$ 
with $\gamma_a = 0.001$,
the steering angle
to $-0.35 \leq \delta \leq 0.35$  
with $\gamma_\delta = 0.001$, 
and the total acceleration 
to 
$0 \leq a_\text{tot}^2 \leq 150$
with $\gamma_{a_\text{tot}^2} = 10^{-8}$,
where $\psi =  25$. The penalizer on the terminal state $X_{K,6}$ 
is zero-mean Gaussian with variance 
$\msgb{V}{X_{K,6}} = 500$. 
The model inputs (\ref{eqn:RaceTrack:DerivativesInput}) 
are unconstrained, which is approximated by 
a zero-mean Gaussian on every $U_k$ with large variance 
$\msgf{V}{U_k} = \diag{10^8, 10^8}$.
The discretization of the race track uses $K=1000$ steps.

\section{Conclusion}

NUP priors allow to incorporate non-Gaussian priors and constraints into linear 
Gaussian models without affecting their computational tractability. 
We proposed new NUP representations 
of half-space constraints,
and we elaborated on recently proposed discretizing NUP priors.
We then discussed the use of such NUP representations
for model predictive control, with a variety of constraints on the input, the 
output, or the internal state of the controlled system. 
In such applications, the computations amount to iterations 
of Kalman-type forward-backward recursions, 
with a complexity (per iteration) that is linear in the planning horizon. 
In consequence, this approach can handle long planning horizons, which 
distinguishes it from the prior art. 
For nonconvex constraints, this approach has no claim to optimality, but it is 
empirically very effective.

The proposed approach was illustrated with a variety of exemplary control problems 
including flappy-bird control and minimal-time race track control. An application 
to a real-world power electronics control problem is demonstrated in a companion 
paper~\cite{keusch2023ContrSysTechn}.

\begin{appendix}
\section{XXXX}

\subsection{Product of Gaussians}

For the convenience of the reader, we state 
\begin{IEEEeqnarray}{rCl}
\IEEEeqnarraymulticol{3}{l}{
\Normal{x; a, \sigma_a^2} \Normal{x; b, \sigma_b^2}
}\nonumber\\\quad
& = & \Normal{x; m_\theta, \sigma^2_\theta} \Normal{a-b; 0, \sigma_a^2 + \sigma_b^2}
      \label{eqn:ProductGaussians}
\end{IEEEeqnarray}
with $\mth$ and $\vth$ as in~(\ref{eqn:box:boxPriorMeanVar}).
For the proof, see 
\cite[Section 1]{bromiley2003products}
or \cite[Appendix~A.1]{Keusch2022PhD}.

\subsection{Proof of Theorem~\ref{theorem:box:hardConstraint}} \label{sec:apds:ProofBox}
\begin{figure}
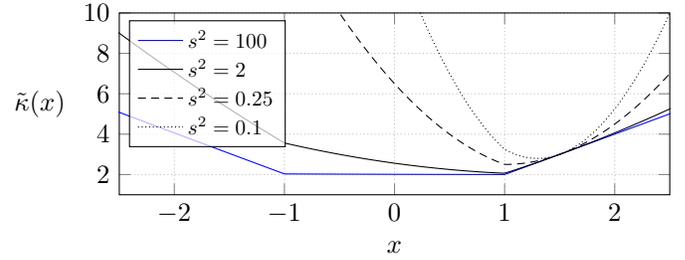

    \centering
    \includestandalone{inc_pdf/costFuncTrivialEx}
    \caption{The cost function (\ref{eqn:box:HardConCostFunc}) 
    for $a=-1, b=1, \mu = 1.5, \gamma=1$, 
    and different values of $s^2$.
    Condition~(\ref{eqn:box:HardCondition}) is satisfied for the solid lines, critically 
    satisfied for the dashed line, and not satisfied for the dotted line.}
    \label{fig:box:costFuncTrivialEx}
  \end{figure}

\noindent
  We first write~(\ref{eqn:box:Analysis:EstX})
  as 
  \begin{IEEEeqnarray}{rCl}
    \hat x 
    &=& \argmax_x p(\breve y \cond x) p_\text{V}(x) \\
    &=& \argmin_x \mleft( - \log\mleft( p(\breve y \cond x)  p_\text{V}(x) \rule{0em}{2ex}\mright) \rule{0em}{2ex}\mright)\\
    & = & \argmin_x \tilde \kappa(x) 
  \end{IEEEeqnarray}
  with
  \begin{equation} \label{eqn:box:HardConCostFunc}
  \tilde\kappa(x) \eqdef \frac{(x-\mu)^2}{2 s^2} + \gamma |x-a| + \gamma |x-b|,
  \end{equation}
  cf.\ Fig.~\ref{fig:box:costFuncTrivialEx}.
  Note that $\tilde \kappa(x)$ is a sum of convex functions and therefore convex itself.
  Consequently, the estimate~(\ref{eqn:box:Analysis:EstX}) 
  is in $[a, b]$ if and only if the global minimum of $\tilde \kappa(x)$ is 
  in $[a, b]$. 
  The latter holds
  if and only if 
  \begin{IEEEeqnarray}{rCl}
    \lim_{\tilde x \uparrow a} \restrict{ \frac{\dd \tilde \kappa(x)}{\dd x}}{x = \tilde x} < 0 
    \quad \text{and} \quad 
    \lim_{\tilde x \downarrow b} \restrict{ \frac{\dd \tilde \kappa(x)}{\dd x}}{x = \tilde x} > 0,
  \end{IEEEeqnarray}
  i.e.,   
  \begin{IEEEeqnarray}{rCl}
    \frac{a-\mu}{s^2} - 2 \gamma < 0 \quad \text{and} \quad \frac{b - \mu}{s^2} + 2 \gamma > 0,
  \end{IEEEeqnarray}
  which boils down to~(\ref{eqn:box:HardCondition}).

\end{appendix}

\balance
\bibliographystyle{IEEEtran}
\bibliography{paper}

\vspace{-2em}

\begin{IEEEbiographynophoto}{Raphael Keusch}
  received the B.Sc. and M.Sc. degrees in electrical engineering from ETH Zurich in 2014 and 2016, respectively. 
  From 2017 to 2018, he was with Sensirion AG, Stäfa, Switzerland. He received the Ph.D. degree in electrical engineering from ETH Zurich in 2022.
  Since 2023, he has been with Verity AG, Zurich, Switzerland.
\end{IEEEbiographynophoto}
  
\vspace{-2em}
\begin{IEEEbiographynophoto}{Hans-Andrea Loeliger}
   received both the Diploma in electrical engineering and the Ph.D. degree (1992) 
  from ETH Zurich, Switzerland. From 1992 to 1995, he was with Linköping University, Linköping Sweden. 
  From 1995 to 2000, he was a technical consultant and coowner of a consulting company. Since 2000, he has 
  been a Professor with the Department of Information Technology and Electrical Engineering of ETH Zurich, 
  Switzerland. His research interests have been in the broad areas of signal processing, machine learning, 
  information theory, communications, error correcting codes, electronic circuits, quantum systems, and 
  neural computation. He is a Fellow of the IEEE.
\end{IEEEbiographynophoto}

\end{document}